\documentclass{amsart} 
\usepackage{a4wide}
\usepackage{amssymb}
\usepackage{graphicx}

\newtheorem{theorem}{Theorem}[section]
\newtheorem{corollary}[theorem]{Corollary}
\newtheorem{lemma}[theorem]{Lemma}

\theoremstyle{definition}

\theoremstyle{remark}

\numberwithin{equation}{section}

\author[E. Duchi]{Enrica Duchi}
\address{E. D.: CAMS, EHESS, 52 bd Raspail, 75006 Paris, France}
\email{Enrica.Duchi@ehess.fr}
\thanks{${}^*$The second author was supported in parts by EC's IHRP
  Program, within the Research Training Network \emph{Algebraic
    Combinatorics in Europe}, grant HPRN-CT-2001-00272.}

\author[G. Schaeffer]{Gilles Schaeffer${}^*$} 
\address{G. S.: LIX,  CNRS, \'Ecole Polytechnique, 91128 Palaiseau, France}
  \email{Gilles.Schaeffer@lix.polytechnique.fr}

\title[Jumping particles]{A combinatorial approach to jumping particles}


\date{October 19, 2003, revised July 25, 2004}

\newcommand{\cross}{\ensuremath{{\!\times\!}}}
\newcommand{\act}{\ensuremath{|\!|}}
\newcommand{\bw}{\ensuremath{{\bullet\atop\circ}}}
\newcommand{\wb}{\ensuremath{{\circ\atop\bullet}}}
\newcommand{\bb}{\ensuremath{{\bullet\atop\bullet}}}
\newcommand{\ww}{\ensuremath{{\circ\atop\circ}}}
\newcommand{\xx}{\ensuremath{{\cross\atop\cross}}}

\begin{document}
\bibliographystyle{plain}

\begin{abstract}
  In this paper we consider a model of particles jumping on a row of
  cells, called in physics the one dimensional totally asymmetric
  exclusion process (TASEP).  More precisely we deal with the TASEP
  with open or periodic boundary conditions and with two or three
  types of particles.  From the point of view of combinatorics a
  remarkable feature of this Markov chain is that it involves Catalan
  numbers in several entries of its stationary distribution.
  
  We give a combinatorial interpretation and a simple proof of these
  observations. In doing this we reveal a second row of cells, which
  is used by particles to travel backward.  As a byproduct we also
  obtain an interpretation of the occurrence of the Brownian excursion
  in the description of the density of particles on a long row of
  cells.
\end{abstract}

\maketitle

\section{Jumping particles}\label{JumpingGilles}

We shall consider a model of jumping particles on a row of $n$ cells
that was exactly solved in the early 90's in physics, under the name
\emph{one dimensional totally asymmetric exclusion process with
open boundaries}, or TASEP for short.
Roughly speaking, the TASEP consists of black particles entering a row
of $n$ cells from an infinite reservoir on the left-hand side and
randomly hopping to the right with the simple exclusion rule that each
cell may contain at most one particle.
\begin{figure}[h]
\begin{center}
\includegraphics[width=70mm]{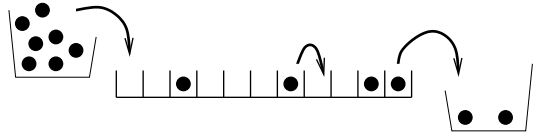}
\end{center}\vspace{-2mm}
\caption{An informal illustration of the TASEP.}
\label{fig:rough}
\end{figure}

The TASEP is usually defined as a continuous-time Markov process on a
finite set of configurations of particles on a line. We shall use an
alternative definition as a finite state Markov chain ---with discrete
time--- which is more convenient for our combinatorial purpose. One
could insist on calling our chain the TASEC, with ``C'' for chain
instead of ``P'' for process, but as we will argue later, there is no
need for this distinction. Another cosmetic modification we allow
ourselves consists in putting a white particle in each empty cell, so
as to make explicit the left-right particle-hole symmetry of the
chain.

\subsection{Definition of the TASEP} 
\begin{figure}[h]
\begin{center}
\includegraphics[width=40mm]{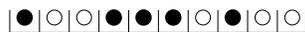}
\end{center}\vspace{-5mm}
\caption{A basic configuration  with $n=10$ cells.}
\label{fig:example-config}
\end{figure}

A \emph{TASEP configuration} is a row of $n$ cells, each containing
either one black particle or one white particle (see
Figure~\ref{fig:example-config}).  These cells are delimited by $n+1$
walls: the left border (or wall $0$), the $i$th separation wall for
$i=1,\ldots,n-1$, and the right border (or wall $n$). 

\begin{figure}
\begin{center}
\includegraphics[width=90mm]{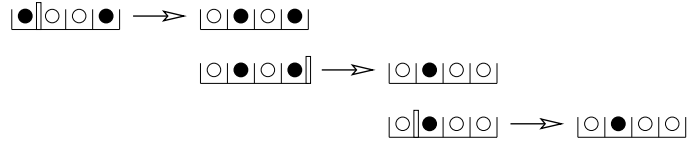}
\end{center}\vspace{-3mm}
\caption{An example of an evolution, with $n=4$ and $\alpha=\beta=\gamma=1$. 
  The active wall triggering each transition is indicated by the
  symbol $\act$.}
\label{fig:example-evolution}
\end{figure}

The TASEP is a Markov chain $S^0_{\alpha\beta\gamma}$ defined on the
set of TASEP configurations for any three parameters $\alpha$, $\beta$
and $\gamma$ in the interval $]0,1]$. From time $t$ to $t+1$, the
chain evolves from the configuration $\tau=S^0_{\alpha\beta\gamma}(t)$
to a configuration $\tau'=S^0_{\alpha\beta\gamma}(t+1)$ as follows:
\begin{itemize}
\item A wall $i$ is chosen uniformly at random among the $n+1$ walls,
  and then may become \emph{active} with probability $\lambda(i)$,
  with $\lambda(i)=\alpha$ for $i=1,\ldots,n-1$, $\lambda(0)=\beta$
  and $\lambda(n)=\gamma$.
\item If the wall does not become active, then nothing happens:
$\tau'=\tau$.
\item Otherwise from $\tau$ to $\tau'$ some changes may occur near the
  active wall:
\begin{itemize}
\item[$a$.] If the active wall is not a border
  ($i\in\{1,\ldots,n-1\}$) and has a black particle on its left-hand
  side and a white one on its right-hand side, then these two particles swap:
  $\bullet\act\circ\to\circ|\bullet$.
\item[$b$.] If the active wall is the left border ($i=0$) and the
  leftmost cell contains a white particle, then this particle
  leaves the system and is replaced by a black one:
  $\act\circ\to|\bullet$.
\item[$c$.] If the active wall is the right border ($i=n$) and the
  rightmost cell contains a black particle, then this particle
  leaves the system and  is replaced by a white one:
  $\bullet\act\to\circ|$.
\item[$d$.] Otherwise the configuration is left unchanged: $\tau'=\tau$.
\end{itemize}
\end{itemize}
As illustrated by Figure~\ref{fig:example-evolution}, black particles
travel from left to right and white particles do the opposite.  The
entire chain for $n=3$ is shown in Figure~\ref{fig:system-basic}.  The
four cases $a,b,c,d$ define an application
$\vartheta:(\tau,i)\mapsto\tau'$ from the set of configurations with
an active wall into the set of configurations. The definition of the
TASEP can be rephrased in terms of this application as: at time $t$
choose a random wall $i=I(t)$ and set
\[
S^0_{\alpha\beta\gamma}(t+1)
=\left\{
\begin{array}{ll}
\vartheta(S^0_{\alpha\beta\gamma}(t),i)  & \textrm{with probability } 
\lambda(i),\\
S^0_{\alpha\beta\gamma}(t) & \textrm{otherwise}.
\end{array}
\right.
\]
The parameters $\alpha$, $\beta$ and $\gamma$ control the rate at
which particles try to move inside the system and at the borders. In
particular, we shall call \emph{maximal flow regime} the special case
$\alpha=\beta=\gamma=1$, in which the rate at which particles try to
move is maximal, and denote $S^0=S^0_{111}$ the corresponding chain.

\subsection{Continuous-time descriptions of the TASEP}
In the physics literature, the TASEP is usually described in the
following terms. The time is continuous, and one considers each wall
independently: during any small time interval $dt$, wall $i$ has probability
$\lambda(i) dt$ to trigger a move $\omega\mapsto\vartheta(\omega,i)$. 
The rate $\lambda(i)$ takes the same values
$\alpha$, $\beta$, $\gamma$ as previously.

Following the probabilistic literature \cite{liggett}, one can give a
formulation which is equivalent to the previous one, but already
closer to ours. In this description, each wall waits for an
independent exponential random time with rate 1 before waking up (in
other terms, at any time, the probability that wall $i$ will still be
sleeping after $t$ seconds is $e^{-t}$). When wall $i$ wakes up, it has
probability $\lambda(i)$ to become active.  If this is the case, then
the transition $\omega\to\vartheta(\omega,i)$ is applied to the
current configuration $\omega$. In any case the wall falls again
asleep, restarting its clock.

This continuous-time TASEP is now easily coupled to the Markov chain
$S^0_{\alpha\beta\gamma}$. Let the time steps of $S^0_{\alpha\beta\gamma}$
correspond to the succession of moments at which a wall wakes up.
Then in both versions, the index of the next wall to wake up is at any
time a uniform random variable on $\{0,\ldots,n\}$, and when a wall
wakes up the transition probabilities are identical. This implies that
the stationary distributions of the continuous-time TASEP and its
Markov chain replica are the same.

\subsection{A remarkable stationary distribution}
Among many results on the TASEP, Derrida \emph{et al.}
\cite{derrida,pasquier} proved the following nice properties of the
chain $S^0=S^0_{111}$, in which particles enter, travel and exit at the
same maximal rate.  First,
\begin{equation}\label{Mafalda}
\textrm{Prob}(S^0(t) \textrm{ contains $0$ black particles})
\;\mathop{\longrightarrow}_{t\rightarrow\infty}\; \frac{1}{C_{n+1}},
\end{equation}
where $C_{n+1}=\frac1{n+2}{2n+2\choose n+1}$ is the $(n+1)$th Catalan
number. More generally, for all $0\leq k\leq n$,
\begin{equation}\label{Sapi}
\textrm{Prob}(S^0(t)\textrm{ contains }k\textrm{ black particles})
\;\mathop{\longrightarrow}_{t\rightarrow\infty}\; 
\frac{\frac{1}{n+1}{n+1\choose k}{n+1\choose n-k}}{C_{n+1}},
\end{equation}
where the numerators are called Narayana numbers.

The model is a finite state Markov chain which is clearly ergodic so
that the previous limits are in fact the probabilities of the same
events in the unique stationary distribution of the chain
\cite{markov}. More generally, Derrida \emph{et al.} provided
expressions for the stationary probabilities of the chain
$S^0_{\alpha\beta\gamma}$ for generic $\alpha$, $\beta$, $\gamma$. Since
their original work a number of papers have appeared providing
alternative proofs and further results on correlations, time
evolutions, etc.  Recent advances and a bibliography can be found for
instance in the article \cite{DLS03}. General books about this kind of
particle processes are \cite{liggett,spohn}.

However, the remarkable appearance of Catalan numbers in the
stationary distribution of $S^0$ is not easily understood from the
proofs in the physics literature. As far as we know, these proofs rely
either on a \emph{matrix ansatz}, or on a \emph{Bethe ansatz}, both
being then proved by a recursion on $n$.

\subsection{Combinatorial results}

Our main ingredient to study the TASEP consists in the construction of
a new Markov chain $X^0_{\alpha\beta\gamma}$ on a set $\Omega^0_n$ of
\emph{complete configurations} that satisfies two main requirements:
on the one hand the stationary distribution of the chain
$S^0_{\alpha\beta\gamma}$ can be simply expressed in terms of that of
the chain $X^0_{\alpha\beta\gamma}$; on the other hand the stationary
behavior of the chain $X^0_{\alpha\beta\gamma}$ is easy to understand.

\begin{figure}[h]
\begin{center}
\begin{minipage}{.4\linewidth}
\includegraphics[width=\linewidth]{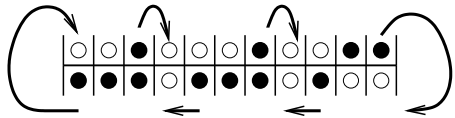}
\end{minipage}
\end{center}
\caption{The circulation of black particles 
  in the complete chain.}
\label{fig:informal}
\end{figure}

The complete configurations that we introduce for this purpose are
made of two rows of $n$ cells containing black and white particles.
The first requirement is met by imposing that, disregarding what
happens in the second row, the chain $X^0_{\alpha\beta\gamma}$
simulates in its first row the chain $S^0_{\alpha\beta\gamma}$.  As
illustrated by Figure~\ref{fig:informal}, the second row will be used
by black and white particles to return to their start point, thus
revealing a circulation of the particles.  The second requirement is
met by adequately choosing the complete configurations and the
transition rules so that $X^0_{\alpha\beta\gamma}$ has a simple
stationary distribution: in particular in the case
$\alpha=\beta=\gamma=1$, $X^0=X^0_{111}$ will have a uniform
stationary distribution.

The chain $X^0_{\alpha\beta\gamma}$ is described in
Section~\ref{sec:completechain}, together with a fundamental property
of its transition rules. Our main result, presented in
Section~\ref{sec:stationarydistribution}, is the combinatorial
interpretation of the stationary distribution of the chain
$S^0_{\alpha\beta\gamma}$, and in particular of
Formulas~(\ref{Mafalda})--(\ref{Sapi}).

\medskip It is known in the literature that some of the
results on the TASEP can be extended to models with three particle
types \cite{ABL,derrida}. We show that this is the case of
Formulas~(\ref{Mafalda})--(\ref{Sapi}) by adapting our approach in
Section~\ref{sec:3tasep} to the 3-TASEP, a Markov chain
$S_{\alpha\beta\gamma\varepsilon}$ in which there are 3 types of
particles, $\bullet$, $\cross$ and $\circ$ and transitions of the form
\[
\bullet\act\cross\to\cross|\bullet, \qquad 
\bullet\act\circ\to\circ|\bullet,
\quad \textrm{and}\quad\cross\act\circ\to\circ|\cross.
\]
Our main results for the 3-TASEP are obtained by a relatively simple
modification of the complete chain. In particular our combinatorial
approach yields the following variant of
Formula~(\ref{Mafalda})--(\ref{Sapi}) for the chain
$S=S_{111\frac12}$: for any $k+\ell+m=n$,
\[
\textrm{Prob}(S(t) \textrm{ contains } \textrm{ $k$ $\bullet$,
$\ell$ $\cross$, and $m$ $\circ$})
 \;\mathop{\longrightarrow}_{t\rightarrow\infty}\;
 \frac{\frac{\ell+1}{n+1}{n+1\choose k}{n+1\choose m}}{\frac{1}{2}
 {2n+2\choose n+1}}.
\]

The TASEP and 3-TASEP are sometimes also defined with periodic
boundary conditions: instead of giving special rules for walls $0$ and
$n$, one identifies these two walls and applies the same rule to every
wall.  With these boundary conditions, the stationary distribution of
the TASEP is easily seen to be uniform.  In Section~\ref{sec:periodic}
we apply our method to study the more interesting distribution of
the 3-TASEP with periodic boundary conditions. In this chain the
number of particles of each type is fixed (since they cannot leave the
system), and, for $k$ $\bullet$, $\ell$ $\times$, $m$ $\circ$, with
$k+\ell+m=n$, we recover the known result:
\[
\textrm{Prob}(\widehat S(t)\;=\;
|\underbrace{\cross\cdots\cross}_\ell
|\underbrace{\circ\cdots\circ}_m
|\underbrace{\bullet\cdots\bullet}_k|)
\;\mathop{\longrightarrow}_{t\rightarrow\infty}\; \frac1{{n \choose
    k}{n \choose m}}.
\]
A different combinatorial proof of this later formula was recently
proposed by Omer Angel \cite{omer}.

\section{The complete chain}\label{sec:completechain}

\subsection{Complete configurations.}
\begin{figure}
\begin{center}
\includegraphics[width=.25\linewidth]{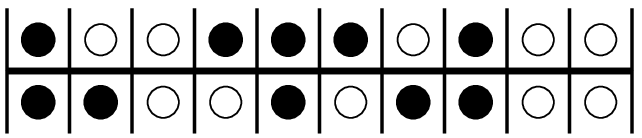}
\hspace{.1\linewidth}
\includegraphics[width=.25\linewidth]{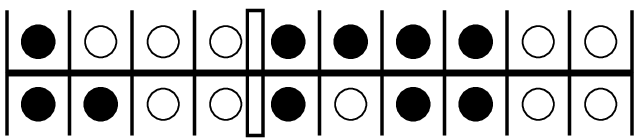}
\end{center}
\vspace{-1mm}
\caption{A complete configuration with $n=10$, and a pair of rows violating the positivity condition at wall 4.}
\label{fig:example-complete-config}
\end{figure}

A \emph{complete configuration} of $\Omega^0_n$ is a pair of rows of
$n$ cells satisfying the following constraints:
\begin{itemize}
\item[\emph{(i)}] \emph{The balance condition}: The two rows contain
  together $n$ black and $n$ white particles.
\item[\emph{(ii)}] \emph{The positivity condition}: On the left of any
  vertical wall there are at least as many black particles as white ones.
\end{itemize}
An example of complete configuration is given in
Figure~\ref{fig:example-complete-config}, together with a pair of rows
that violates the positivity condition.

Given a complete configuration of length $n$, and an integer $j$,
$0\leq j\leq n$, let $B(j)$ and $W(j)$ be respectively the numbers of
black and white particles lying in the first $j$th columns (from left
to right), and set $E(j)=$ $B(j)-W(j)$. In other terms, the quantities
$B(j)$, $W(j)$ and $E(j)$ represent the number of black particles, the
number of white particles, and their difference on the left
of wall $j$. In particular, $E(0)=E(n)=0$, and Condition~($ii$)
of the definition of complete configurations reads $E(j)\geq0$ for
$j=0,\ldots,n$ (this is why we call it a positivity condition).
Readers with background in enumerative combinatorics may have recognized
here complete configurations with $n$ columns as bicolored Motzkin
paths with $n$ steps, or Dyck paths with $2n+2$ steps in disguise
\cite[Chap.  6]{stanley}. 
In particular these characterizations yield
the following lemmas. A direct proof of these lemmas is given in
Section~\ref{sec:counting} for completeness.
\begin{lemma}\label{lemma:olaf}
  The number $|\Omega^0_{n}|$ of complete configurations is
  $C_{n+1}=\frac{1}{n+1} { 2n+2 \choose n}=\frac1{n+2}{2n+2\choose n+1}$.
\end{lemma}
\begin{lemma}\label{lemma:ciclico}
  Let $k,m,n$ be non negative integers with $k+m=n$.  The
  number $|\Omega^0_{k,m}|$ of complete configurations of
  $\Omega^0_{n}$ with $k$ black and $m$
  white particles on the top row, and $m$ black and $k$ white
  particles on the bottom row is $\frac{1}{n+1}{n+1\choose
    k}{n+1\choose m}$.
\end{lemma}
A first hint to our interest in complete configurations should follow
from the comparison of the lemmas with the probabilities
in~(\ref{Mafalda}) and~(\ref{Sapi}).

\subsection{First definition of the complete chain}\label{sec:informal}
\begin{figure}
  \begin{center}
    \includegraphics[width=90mm]{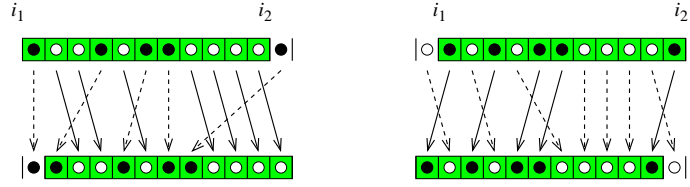}
    \vspace{-3mm}
  \end{center}
  \caption{A white sweep and a black sweep.}
  \label{fig:interpretation}
    \vspace{-1mm}
\end{figure}

The Markov chain $X^0_{\alpha\beta\gamma}$ on $\Omega^0_n$ will be
defined in terms of an application $T$ from the set $\Omega^0_n \times
\{0,\ldots,n\}$ to the set $\Omega^0_n$ that extends the application
$\vartheta$.
Given a complete configuration $\omega$ and an active wall $i$, the
action of $T$ on the top row of $\omega$ does not depend on the second
row, and mimics the action of $\vartheta$ as defined by cases $a$,
$b$, $c$ and $d$ of the description of the TASEP.  In particular in
the top row, black particles travel from left to right and white
particles from right to left.  As opposed to that, in the bottom row,
$T$ moves black and white particles backward, as illustrated by
Figure~\ref{fig:informal}.
In order to describe how these moves are performed, we first introduce the
concept of sweep (see Figure~\ref{fig:interpretation}):
\begin{itemize}
\item A \emph{white sweep} between walls $i_1$ and $i_2$ consists in
  all white particles that are in the bottom row between walls $i_1$
  and $i_2$ simultaneously hopping to the right (some black particles
  thus being displaced to the left in order to fill the gaps). For
  well definiteness a white sweep between $i_1$ and $i_2$ can occur
  only if the particle on the right-hand side of $i_2$ is black.
\item A \emph{black sweep} between walls $i_1$ and $i_2$ consists in
  all black particles that are in the bottom row between walls $i_1$ and
  $i_2$ simultaneously hopping to the left (some white particles thus
  being displaced to the right in order to fill the gaps). For well
  definiteness a white sweep between $i_1$ and $i_2$ can occur only if
  the particle on the left-hand side of $i_1$ is white.
\end{itemize}
Next, given a complete configuration and a wall $i$, we distinguish
the following walls: if there is a black particle on the left-hand
side of wall $i$ in the top row, let $j_1<i$ be the leftmost wall such that
there are only white particles in the top row between walls $j_1$ and
$i-1$; if there is a white particle on the right-hand side of $i$ in
the top row, let $j_2>i$ be the rightmost wall such that there are
only black particles in the top row between walls $i+1$ and $j_2$.

\begin{figure}
\begin{center}
  \includegraphics[width=110mm]{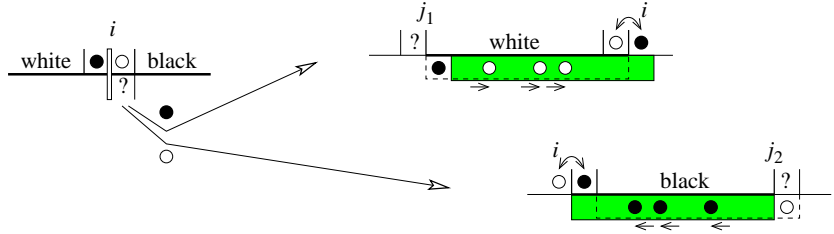}
\end{center}
\vspace{-4mm}
\caption{Sweeps occurring below the 
  transition $(\bullet\act\circ\to\circ|\bullet)$.}
\label{fig:ambiguous}
\end{figure}
With these definitions, we are in the position to describe the action
of $T$.  Given a complete configuration $\omega\in\Omega^0_n$ and a wall
$i\in\{0,\ldots,n\}$, the cases $a$, $b$, $c$ and $d$ of the
transition rule $\vartheta$ describe the top row of the image
$T(\omega,i)$, and they are complemented as follows to describe the
bottom row of the image:
\begin{itemize}
\item[$a$.] In this case $i\in\{1,\ldots,n-1\}$ and this wall
  separates a black and a white particle in the top row of $\omega$.
  The moves in the bottom row then depend on the particle on the
  bottom right of wall $i$ in $\omega$: if it is black, a white sweep
  occurs between $j_1$ and $i$, otherwise it is white and a black
  sweep occurs between $i+1$ and $j_2+1$ (or between $i+1$ and $n$
  if $j_2=n$).  These moves are illustrated by
  Figure~\ref{fig:ambiguous} (see also
  Figures~\ref{move-right-sweep}--\ref{move-left-sweep}, left and middle).
\item[$b$.] In this case $i=0$ and the leftmost particle of the top
  row of $\omega$ is white. Then the leftmost column of $\omega$ is a
  $|\wb|$-column. These two particles exchange (so that a black
  particle enters in the top row in agreement with rule $b$ for
  $\vartheta$), and a black sweep occurs between the
  left border and wall $j_2+1$, or between the left and right borders
  if $j_2=n$ (see Figure~\ref{move-left-sweep}, right).
\item[$c$.] In this case $i=n$ and the rightmost particle of the top
  row of $\omega$ is black. Then the rightmost column of $\omega$ is a
  $|\bw|$-column. These two particles exchange (so that a white
  particle enters in the top row in agreement with rule $c$ for
  $\vartheta$), and a white sweep occurs between wall $j_1$ and the
  right border (see Figure~\ref{move-right-sweep}, right).
\item[$d$.] Otherwise nothing happens.
\end{itemize}
The fact that the configuration $T(\omega,i)$ produced in each case
satisfies the positivity constraint is not difficult to prove and it
is explicitly checked in the next section.

The Markov chain $X^0_{\alpha\beta\gamma}$ on the set $\Omega^0_n$ of
complete configurations with length $n$ is defined from $T$ exactly as
the TASEP is described from $\vartheta$: the evolution rule from time
$t$ to $t+1$ consists in choosing $i=I(t)$ uniformly at random in
$\{0,\ldots,n\}$ and setting
\[
X^0(t+1)=\left\{
\begin{array}{ll}
T\left(X^0(t),i\right) &\textrm{ with probability $\lambda(i)$,}\\
X^0(t) &\textrm{ otherwise.}
\end{array}\right.
\]

By construction, the Markov chains $S^0_{\alpha\beta\gamma}$ and
$X^0_{\alpha\beta\gamma}$ are related by
\[
S^0_{\alpha\beta\gamma}\;\equiv\;\textrm{top}(X^0_{\alpha\beta\gamma}),
\]
where $\textrm{top}(\omega)$ denotes the top row of a complete
configuration $\omega$, and the $\equiv$ is intended as identity in
law at any time, provided $S_{\alpha\beta\gamma}^0(0)$ and
$\textrm{top}(X_{\alpha\beta\gamma}^0(0))$ are equally distributed. 

\begin{figure}
\begin{center}
\begin{minipage}{.6\linewidth}
  \includegraphics[width=\linewidth]{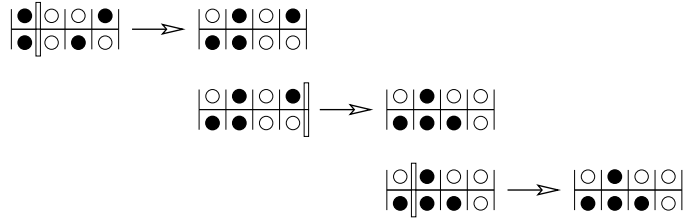}
\end{minipage}
\end{center}
\vspace{-2mm}
\caption{An  example of actual evolution with $n=4$ 
and $\alpha=\beta=\gamma=1$.}
\label{fig:example-evolution-comp}
\end{figure}

An appealing interpretation from a combinatorial point of view is that
we have revealed a circulation of the particles, that use the bottom
row to travel backward and implement the infinite reservoirs, as
illustrated by Figure~\ref{fig:informal}.  An example of evolution is
given by Figure~\ref{fig:example-evolution-comp}. The TASEP and
complete chain with two particles for $n=3$ are represented in
Figures~\ref{fig:system-basic}--\ref{fig:system-complete}.

\subsection{Restatement of the transition rules: the bijection $\bar T$}
\label{sec:bijection-intro}
\begin{theorem}\label{thm:bijection}
  The application $T$ is the first component
  $\Omega^0_n\times\{0,\ldots,n\}\to\Omega^0_n$ of a bijection $\bar
  T$ from $\Omega^0_n\times\{0,\ldots,n\}$ into itself.
\end{theorem} 
\begin{proof}
  In order to define the application $\bar T$, we shall partition the
  set $\Omega^0_n\times\{0,\ldots,n\}$ into classes $A_{a'}$,
  $A_{a''}$, $A_{b}$, $A_{c}$, and $A_{d}$, and describe, for each
  class $A$, its image $B=\bar T(A)$. From
  now on in this section, $(\omega,i)$ and $(\omega',j)$ respectively
  denote an element of the current class and its image, and $j_1$ and
  $j_2$ are defined from $(\omega,i)$ as in
  Section~\ref{sec:informal}. 
  
  We are going to describe the image $(\omega',j)$ of $(\omega,i)$ by
  $\bar T$ in terms of deletions and insertions of $|\bw|$- or
  $|\wb|$-columns or of ${\bullet\atop}|{\atop\circ}$-diagonals. One
  advantage of these operations is that they clearly preserve the
  balance and positivity conditions, so we will directly know in each
  case that the image $\omega'$ belongs $\Omega^0_n$.  The reader is
  invited to chek, using Figures~\ref{move-right-sweep}
  and~\ref{move-left-sweep}, that the configuration $\omega'$ obtained
  in each case is, as claimed in the theorem, the same as the
  configuration $T(\omega,i)$ that was described in terms of sweeps in
  the previous section:
\begin{itemize} 
\item %[$A_{a_1}$] 
  \textsf{If the wall $i$ separates in the top row of $\omega$ a black
    particle $P$ and a white particle $Q$.}  There are two cases
  depending on the type of the particle $R$ that is below $Q$ in
  $\omega$:
  \begin{itemize}
  \item[$A_{a'}$] \textsf{The particle $R$ is black.}  Then $j=j_1$
    and $\omega'$ is obtained by moving the $|\wb|$-column $|{Q\atop
      R}|$ from the right-hand side of wall $i$ to the right-hand side
    of wall $j$ (Figure~\ref{move-right-sweep}, left-middle).
    
    The image $\smash{B_{a'}}$ of the class $A_{a'}$ consists of pairs
    $(\omega',j)$ such that: the wall $j$ is the left border ($j=0$)
    or it has a black particle on its left-hand side in the top row,
    there is a $|\wb|$-column on the right-hand side of wall $j$, and
    the sequence of white particles on the right-hand side of wall $j$
    in the top row is followed by a black particle.
    
  \item[$A_{a''}$] \textsf{The particle $R$ is white}.  Then $j=j_2$
    and $\omega'$ is obtained by moving the particles $P$ and $R$ from
    wall $i$ (where they form a
    ${\bullet\atop}|{\atop\circ}$-diagonal) to wall $j$ so that they
    form a ${\bullet\atop}|{\atop\circ}$-diagonal if $j<n$
    (Figure~\ref{move-left-sweep}, left), or a $|\bw|$-column if $j=n$
    (Figure~\ref{move-left-sweep}, middle).

    The image $B_{a''}$ of the class $A_{a''}$ consists
    of pairs $(\omega',j)$ with a $|\ww|$-column
    or the border on the right-hand side of wall $j$ of
    $\omega'$ and such that there is a non-empty sequence of black
    particles on the left-hand side of wall $j$ in the top row,
    followed by a white particle.
  \end{itemize}

  \begin{figure}
    \begin{center}
     \includegraphics[width=\linewidth]{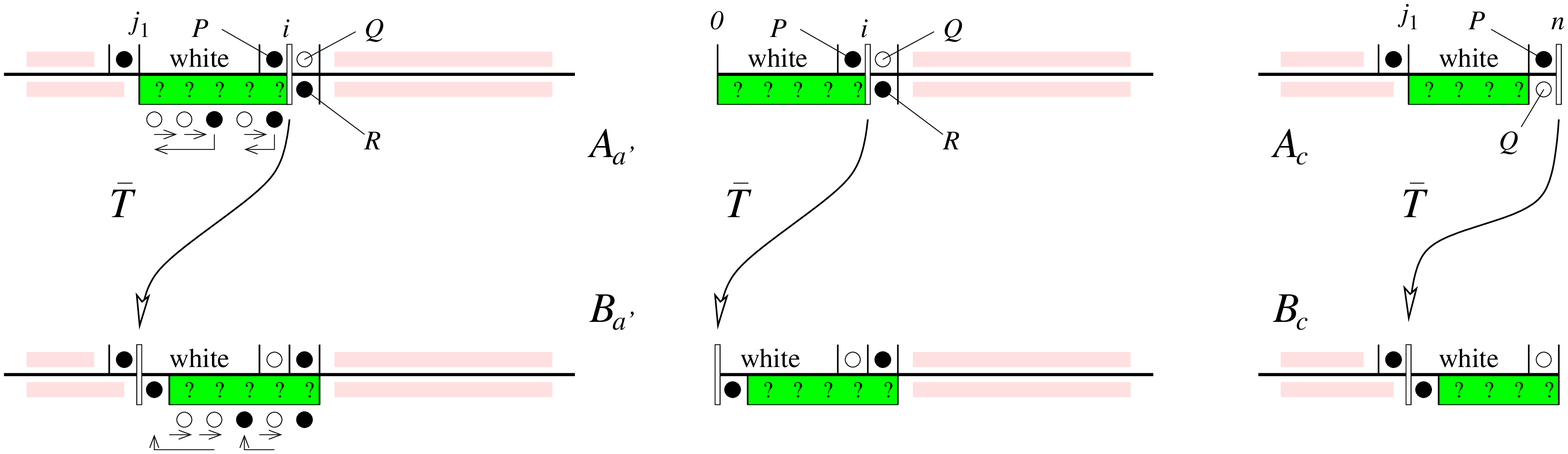}
    \end{center}
    \caption{Moves in cases $A_{a'}$ and $A_{c}$. Below the two
      left-hand side configurations, the white sweep in the
      bottom row is illustrated on an exemple.}
    \label{move-right-sweep}
    \begin{center}
     \includegraphics[width=\linewidth]{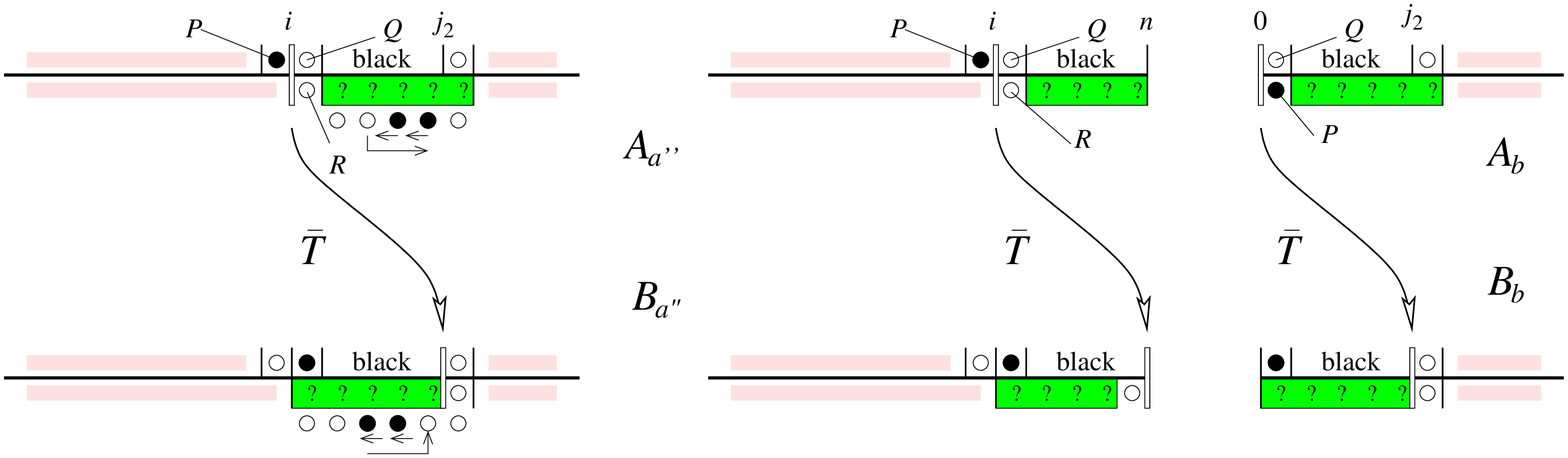}
    \end{center}
    \caption{Moves in the cases $A_a''$ and $A_{b}$. Below the two 
      left-hand side configurations, the black sweep in the bottom row 
      is illustrated on an exemple.}
    \label{move-left-sweep}
  \end{figure} 
    
\item[$A_{b}$] \textsf{If $i=0$ and there is a white particle $Q$ in
    the leftmost top cell of $\omega$.}  The cell under $Q$ then
  contains a black particle $P$ (Figure~\ref{move-left-sweep}, right).
  Then $j=j_2$ and $\omega'$ is obtained by moving $P$ and $Q$ to wall
  $j$ so that they form a ${\bullet\atop}|{\atop\circ}$-diagonal if
  $j<n$ or a $|\bw|$-column if $j=n$.
  
  The image $B_{b}$ of the class $A_{b}$ consists of pairs
  $(\omega',j)$ with a $|\ww|$-column or the border on the right-hand side of
  wall $j$ of $\omega'$ and such that there is a non-empty
  sequence of black particles on the left-hand side of wall $j$ in the
  top row, ending at the left border.

%%%%%%%%%
\item[$A_{c}$] \textsf{If $i=n$ and there is a black particle $P$ in
    the rightmost top cell of $\omega$.}  The cell under $P$ then
  contains a white particle $Q$ (Figure~\ref{move-right-sweep},
  right). Then $j=j_1$ and $\omega$ is obtained by moving $P$ and $Q$ to 
  wall $j$ so that they form a $|\wb|$-column on its right-hand side.    
  
  The image $B_{c}$ of the class $A_{c}$ consists of pairs
  $(\omega',j)$ such that: the wall $j$ is the left border ($j=0$) or
  it has a black particle on its left-hand side in the top row, there
  is a $|\wb|$-column on the right-hand side of wall $j$, and the
  sequence of white particles on the right-hand side of wall $j$ in
  the top row ends at the right border.
%%%%%

\item[$A_d$] \textsf{This class contains all the remaining pairs
    $(\omega,i)$}.  These configurations are left unchanged by $\bar
  T$, so that $B_d=\bar T(A_d)=A_d$.
\end{itemize}
In each case of the definition, the transformation described is
reversible: from $(\omega',k)$ in one of the image classes, the wall
$i$ and then the configuration $\omega$ are easily recovered. The
theorem thus follows from the fact that $\{B_{a'},B_{a''},
B_{b},B_{c}, B_{d}\}$ is a partition of
$\Omega^0_n\times\{0,\ldots,n\}$.
\end{proof}

\section{Stationary distributions}\label{sec:stationarydistribution}

The Markov chain $X^0_{\alpha\beta\gamma}$ is clearly aperiodic and we
check in Section~\ref{sec:paths} that it is irreducible, \emph{i.e.}
that there is an evolution between any two configurations.  This
implies that the chain $X^0_{\alpha\beta\gamma}$ is ergodic, \emph{i.e.}
it has a unique stationary distribution, to which
$X^0_{\alpha\beta\gamma}(t)$ converges as $t$ goes to infinity \cite{markov}.
Our aim in this section is to find this distribution and to use it to
give a combinatorial interpretation to that of $S^0_{\alpha\beta\gamma}$.

We first deal with the maximal flow regime, for which all ingredients
are now ready. Then we discuss the generic case.

\subsection{The maximal flow regime $\alpha=\beta=\gamma=1$}
\label{sec:maximal}
\begin{theorem}\label{thm:proba}
The Markov chain $X^0$ has a uniform stationary distribution.
\end{theorem}
\begin{figure}
%\vspace{3em}
\begin{center}
   \includegraphics[width=.8\linewidth]{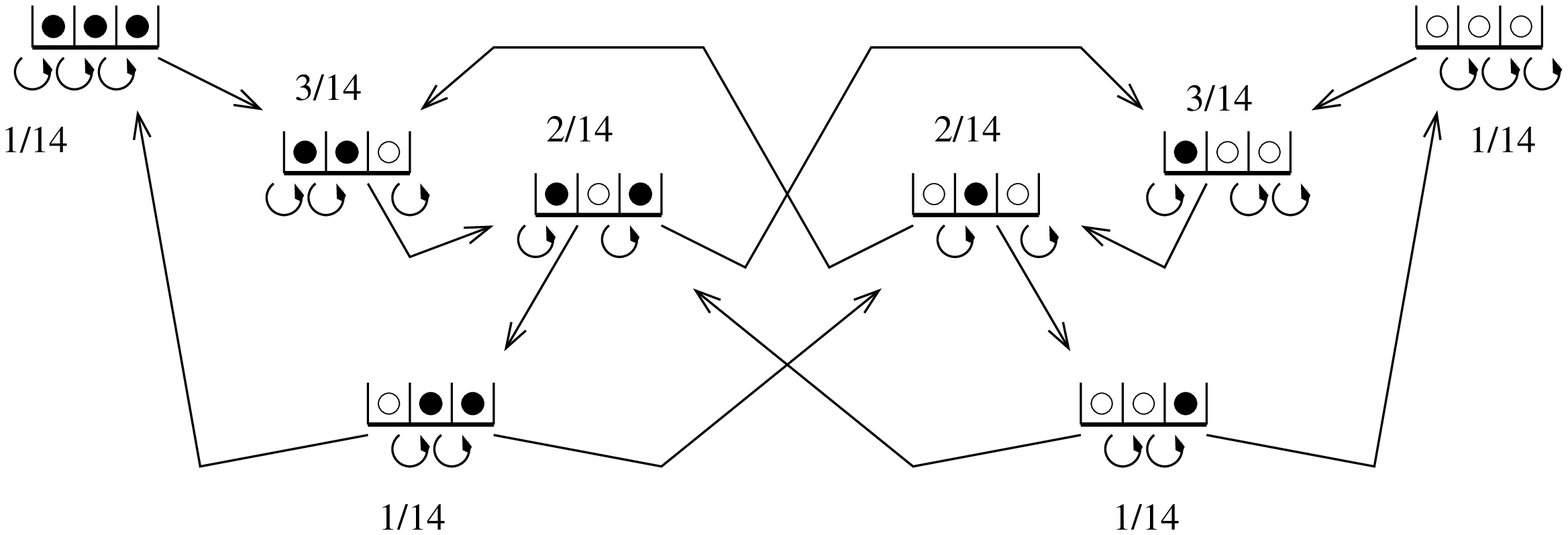}
 \end{center}
 \caption{The TASEP configurations for $n=3$ and transitions between
   them. The starting point of each arrow indicates the wall
   triggering the transition.  The numbers are the stationary
   probabilities for $\alpha=\beta=\gamma=1$.}
 \label{fig:system-basic}
\vspace{1em}
 \begin{center}
 \includegraphics[width=.9\linewidth]{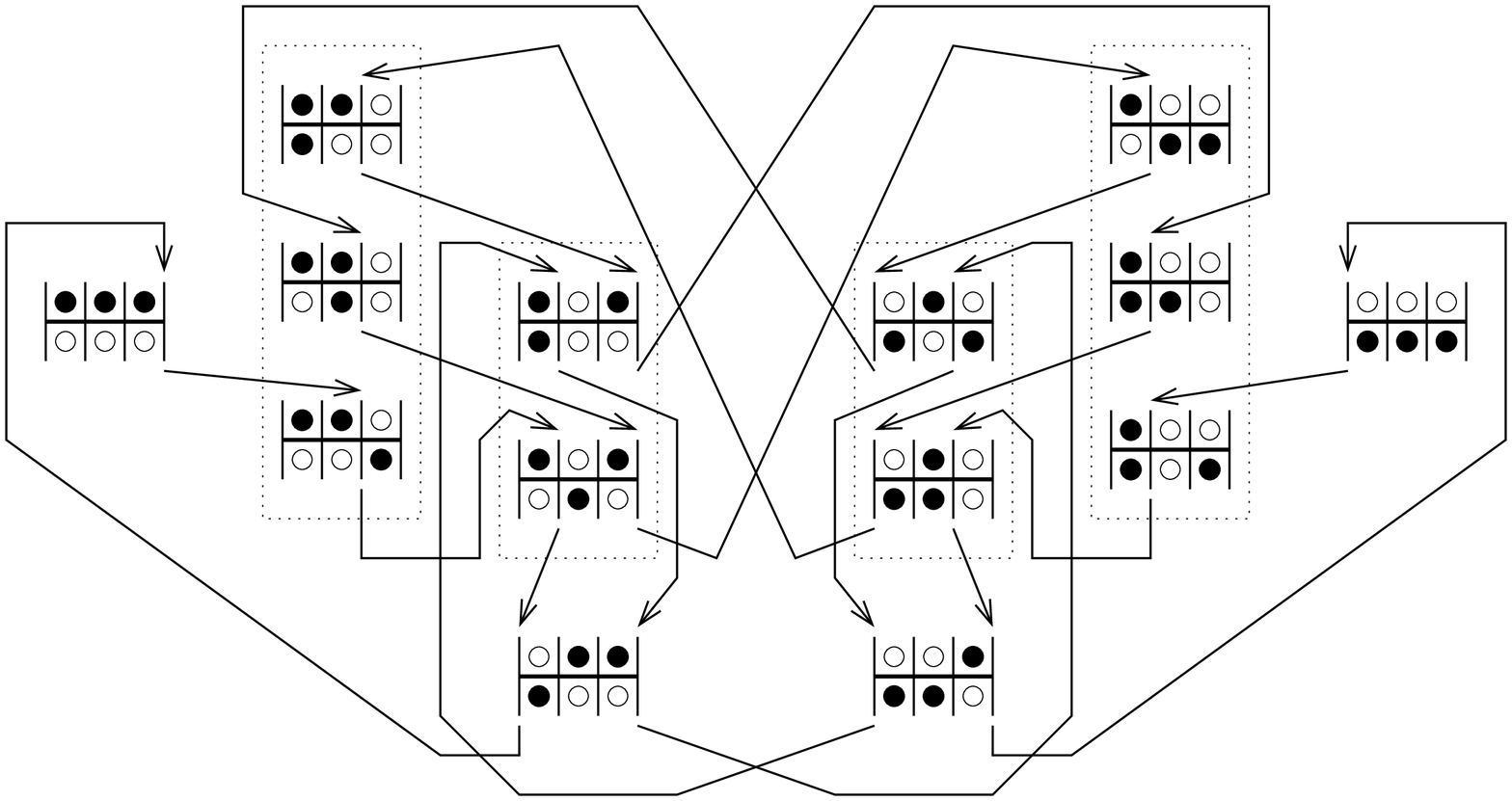}
 \caption{The $14$ complete configurations for $n=3$ and transitions
   between them. The starting point of each arrow indicates the wall
   triggering the transition (loop transitions are not indicated).
   For $\alpha=\beta=\gamma=1$, the stationary probabilities are
   uniform (equal to $1/14$) since each configuration has equal in-
   and out-degrees. Ignoring the bottom rows reduces this Markov chain
   to the chain of Figure~\ref{fig:system-basic}.}
 \label{fig:system-complete}
 \end{center}
%\vspace{3em}
\end{figure}
\begin{proof}
  As illustrated by Figure~\ref{fig:system-complete},
  Theorem~\ref{thm:bijection} says that the vertices of the
  transition graph of the chain have equal in- and
  out-degrees. Moreover the $n+1$ possible transitions from a
  configuration $\omega$ have equal probabilities to be chosen, since
  the active wall is chosen uniformly in $\{0,\ldots,n\}$. The uniform
  distribution on $\Omega^0_n$ hence clearly satisfies the local
  stationarity equation at each configuration $\omega$: assuming that
  at time $t$ the distribution is uniform,
\[
\textrm{Prob}(X(t)=\omega)=\frac1{|\Omega^0_n|}\qquad\textrm{for all}\;\omega,
\]
then at time $t+1$, it remains uniform, since
\[
{\textrm{Prob}(X(t+1)=\omega')}
=\!\!\!
\sum_{(\omega,i)\in  T^{-1}(\omega')} \!\!\!\!
\textrm{Prob}(X(t)=\omega)\cdot{\frac1{n+1}}
\;=\;
{\big| T^{-1}(\omega')\big|}\cdot\frac1{|\Omega^0_n|}\cdot\frac{1}{n+1}
\;=\;\frac1{|\Omega^0_n|},
\]
where $T^{-1}(\omega')$ denotes the set of preimages of $\omega'$
respectively by $T$. The last equality follows from the facts that
$T^{-1}(\omega')= \{ \bar T^{-1}(\omega', j) \mid j=0, \ldots, n\}$,
and that $\bar T$ is a bijection.
\end{proof}

The relation $S^0\equiv\textrm{top}(X^0)$ now allows us to derive from
Theorem~\ref{thm:proba} the annonced combinatorial interpretation of
Formulas~(\ref{Mafalda}) and~(\ref{Sapi}).
\begin{theorem}
  Let $\textrm{top}(\omega)$ denote the top row of a complete
  configuration $\omega$. Then for any initial distribution $S^0(0)$
  and $X^0(0)$ with $S^0(0)\equiv\textrm{top}(X^0(0))$, and any TASEP
  configuration $\tau$,
 \[
 \textrm{Prob}(S^0(t) = \tau)\;=\;
 \textrm{Prob}(\textrm{top}(X^0(t)) =\tau)
 \;\mathop{\longrightarrow}_{t\rightarrow\infty}\;
 \frac{\big|\{\omega\in\Omega^0_n\mid \textrm{top}(\omega)=\tau\}|}
 {|\Omega^0_{n}|}.
 \]
In particular, for any $k+m=n$, we obtain combinatorially the formula:
\[
\textrm{Prob}(S^0(t) \textrm{ contains } k \textrm{ black and }
 m \textrm{ white particles})
\;\mathop{\longrightarrow}_{t\rightarrow\infty}\;
 \frac{|\Omega^0_{k,m}|}{|\Omega^0_n|}
\;=\;
 \frac{\frac{1}{n+1}{n+1\choose k}{n+1\choose m}}{C_{n+1}}.
 \]
\end{theorem}
As discussed in Section~\ref{sec:conclusion}, this interpretation sheds
a new light on some recent results of Derrida \emph{et al.} connecting
the TASEP to Brownian excursions \cite{Derrida03brownian}.

\subsection{Arbitrary $\alpha$, $\beta$ and $\gamma$}
\label{sec:arbitrary}

In order to express the stationary distribution of the general chain
$X^0_{\alpha\beta\gamma}$, we associate a weight $q(\omega)$ to each
complete configuration $\omega$, which is defined in terms of two
combinatorial statistics.

By definition, a complete configuration $\omega$ is a concatenation of
four types of columns $|{\bullet\atop\bullet}|$,
$|{\bullet\atop\circ}|$, $|{\circ\atop\bullet}|$ and
$|{\circ\atop\circ}|$, subject to the balance and positivity
conditions.  In particular, the concatenation of two complete
configurations of $\Omega^0_i$ and $\Omega^0_j$ with $i+j=n$ yields a
complete configuration of $\Omega^0_n$. Let us call \emph{prime} a
configuration that cannot be decomposed in this way. A complete
configuration $\omega$ can be uniquely written as a concatenation
$\omega=\omega_1\cdots\omega_m$ of prime configurations. These prime
factors can be of three types: $|\bw|$-columns, $|\wb|$-columns, and
\emph{blocks} of the form $|\bb|\omega'|\ww|$ with $\omega'$ a
complete configuration.  The inner part $\omega'$ of a block
$\omega=|\bb|\omega'|\ww|$ is referred to as its \emph{inside}.

\begin{figure}
\begin{center}
 \includegraphics[width=70mm]{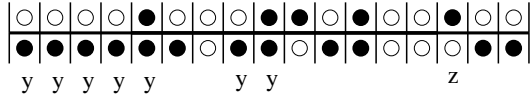}
\end{center}
\vspace{-4mm}
\caption{A configuration $\omega$ with weight $q(\omega)=\alpha^8\beta^{10}\gamma^{16}$. Labels are indicated below particles. }
\label{fig:qpar}
\end{figure}

Now, given a complete configuration $\omega$, let us assign labels to
some of the particles of its bottom row: first, each white particle is
labeled $z$ if it is not in a block, and then, each black particle
is labeled $y$ if it is not in the inside of a block and there are no
$z$ labels on its left.  The number of labels of type $y$ and the
number of labels of type $z$ in a configuration $\omega$ will be
denoted $n_y(\omega)$ and $n_z(\omega)$ respectively. Then the
\emph{weight} of a configuration $\omega$ is defined as
\[
q(\omega)\;=\;{\beta^n\gamma^n}
\left(\frac{\alpha}{\beta}\right)^{n_y(\omega)}
\left(\frac{\alpha}{\gamma}\right)^{n_z(\omega)}
\;=\;\alpha^{n_y(\omega)+n_z(\omega)}
\beta^{n-n_y(\omega)}\gamma^{n-n_z(\omega)}.
\]
In other terms, there is a factor $\alpha$ per label, a factor $\beta$
per unlabeled black particle and a factor $\gamma$ per unlabeled white
particle.  For instance, the weight of the configuration of
Figure~\ref{fig:qpar} is $\alpha^{8}\beta^{10}\gamma^{16}$, and more
generally the weight is a monomial with total degree $2n$.

\begin{theorem}\label{thm:generic}
  The Markov chain $X^0_{\alpha\beta\gamma}$ has the following unique
  stationary distribution:
  \[
  \textrm{Prob}(X^0_{\alpha\beta\gamma}(t)=\omega)
  \;\mathop{\longrightarrow}_{t\rightarrow\infty}\;
  \frac{q(\omega)}{Z^0_n}, \qquad
  \textrm{where } Z^0_n={\sum_{\omega' \in \Omega^0_n}q(\omega')},
  \]
  where $q(\omega)$ is the previously defined weight on complete
  configurations.
\end{theorem}
Since $X^0_{\alpha\beta\gamma}$ is aperiodic and irreducible, it is
sufficient to show that the distribution induced by the weights $q$ is
stationary.  The result is based on a further property of the
bijection $\bar T$ of which $T$ is the first component. 
\begin{lemma} \label{lem:preserve}
  The bijection $\bar
    T:\Omega^0_n\times\{0,\ldots,n\}\to\Omega^0_n\times\{0,\ldots,n\}$
    transports the weights:
    \begin{equation}\label{for:preserve}
    \lambda(i)\,q(\omega)\;=\;\lambda(j)\,q(\omega'),\qquad
    \textrm{for all $(\omega',j)=\bar T(\omega,i)$,}
    \end{equation}
    where $\lambda(i)=\alpha$ for $i\in\{1,\ldots,n-1\}$,
    $\lambda(0)=\beta$ and $\lambda(n)=\gamma$.
  \end{lemma}
\begin{proof}

Let $\omega$ be a complete configuration belonging to $\Omega_n^0$. The
following properties will be useful:
\begin{itemize}
\item \textsf{Property 1.}  \emph{In a local configuration
  $|{\bullet\atop?}|\wb|$ the black particle in the bottom row never
  contributes a label $y$.}

The black particle of a $|\wb|$-column can contribute a label $y$ only
if it is not in the inside of a block. This happens only if the
particle $?$ is white and is not in a block. But then this particle
carries a label $z$ which is on the left of the black particle.

\item \textsf{Property 2.} \emph{The bottom white particle of a
  $|\ww|$-column never contributes a label $z$.}

 This property is immediate since a
  $|\ww|$-column is always in a block.

\item \textsf{Property i.} \emph{The deletion/insertion of a $|\wb|$-column 
  does not change the labels of other particles.}

When a $|\wb|$-column is inserted or removed in the inside of a block,
the block structure is unchanged and there is no effect on labels.
When it is inserted or removed at a position no included in a block it
may contribute a label $y$, but this has no effect on other labels.

\item \textsf{Property ii.}\emph{ The deletion/insertion of a
  ${\bullet\atop}|{\atop\circ}$-diagonal taking the form
  $|{\bullet\atop?}|\ww|\leftrightarrow|{\circ\atop?}|$ does not
  change the labels of other particles.}

 The situation $|\bw|\ww|\leftrightarrow |\ww|$
  can be view as the insertion or deletion of a $|\bw|$-column in the inside
  of a block, which has no effect. The other situation
  $|\bb|\ww|\leftrightarrow|\wb|$ may occur outside a block, in which
  case a white particle is added or removed on the bottom row, but
  in a small block $|\bb|\ww|$.
\end{itemize}

  The relation is checked using these properties by comparing
  $q(\omega)$ and $q(\omega')$ in each case of the definition of the
  bijection $\bar T$.
  \begin{itemize}
  \item[$A_{a'}$] If $j\neq0$ (Figure~\ref{move-right-sweep}, left),
    according to Property~1 the particle $R$ does not contribute a
    label $y$ neither in $\omega$ nor in $\omega'$.  Moreover,
    according to Property~i, the displacement of the $|\wb|$-column
    does not affect labels of other particles. Hence
    $q(\omega)=q(\omega')$, in agreement with $\lambda(i)=\lambda(j)=\alpha$.
    
    If instead $j=0$ (Figure~\ref{move-right-sweep}, middle),
    Property~1 applies only to $\omega$: in $\omega'$, the displaced
    $|\wb|$-column is the leftmost one, so that its black particle
    contributes a supplementary $y$ label.  Therefore
    $q(\omega')=q(\omega)\frac{\alpha}{\beta}$, in agreement with
    $\lambda(i)=\alpha$, $\lambda(0)=\beta$.
      
  \item[$A_{a''}$] If $j\neq n$ (Figure~\ref{move-left-sweep}, left),
      from Property~2 we see that the particle $R$ does not contribute
      a label $y$ neither in $\omega$ nor in $\omega'$. Observe
      moreover that the displacement of a
      ${\bullet\atop}|{\atop\circ}$-diagonal does not affect
      labels of other particles according to Property ii. Hence
      $q(\omega)=q(\omega')$, in agreement with
      $\lambda(i)=\lambda(j)=\alpha$.
      
      If $j=n$ (Figure~\ref{move-left-sweep}, middle), Property~2
      applies only to $\omega$: the move amounts to deleting a
      ${\bullet\atop}|{\atop\circ}$-diagonal and inserting a
      $|\bw|$-column at the right border. The white particle of this column
      thus contributes a $z$ label. Therefore
      $q(\omega')=q(\omega)\frac{\alpha}{\gamma}$, in agreement with
      $\lambda(i)=\alpha$ and $\lambda(n)=\gamma$.

  \item[$A_{b}$] If $j\neq n$ (Figure~\ref{move-left-sweep}, right),
    the move consists in deleting a $|\wb|$-column, which is the
    leftmost and thus contributes a $y$ label in $\omega$, and
    inserting a ${\bullet\atop}|{\atop\circ}$-diagonal, which
    according to Property~2 does not contribute a $z$ label. According
    to Property~i and~ii the other labels are left unchanged.
    Therefore $q(\omega')=q(\omega)\frac{\beta}{\alpha}$, in agreement
    with $\lambda(0)=\beta$ and $\lambda(j)=\alpha$.

    If $j=n$, $q(\omega')=q(\omega)\frac{\beta}{\alpha}
    \frac{\alpha}{\gamma}= q(\omega)\frac{\beta}{\gamma}$, in
    agreement with $\lambda(0)=\beta$ and $\lambda(n)=\gamma$.
    
  \item[$A_c$] If $j\neq 0$ (Figure~\ref{move-right-sweep}, right),
    $\omega'$ is obtained by deleting a $|\bw|$-column on the
    left-hand side of the wall $n$ and inserting a $|\wb|$-column on
    the right-hand side of $j_1$.  According to Property~i only the labels of
    displaced particles can be affected. Since the deleted
    $|\bw|$-column is the rightmost column, its white particle
    contributes a $z$ label in $\omega$.  As opposed to that,
    Property~1 forbids the $|\wb|$-column to contribute a label in
    $\omega'$.  Therefore $q(\omega')=q(\omega)\frac{\gamma}{\alpha}$,
    in agreement with $\lambda(n)=\gamma$ and $\lambda(j)=\alpha$.

    If $j=0$,  $q(\omega')=q(\omega)\frac{\gamma}{\alpha}
    \frac{\alpha}{\beta}= q(\omega)\frac{\gamma}{\beta}$, in agreement
    with $\lambda(n)=\gamma$ and $\lambda(0)=\beta$.
\end{itemize}
\end{proof}
\begin{proof}[Proof of Theorem~\ref{thm:generic}.]
  In order to see that the distribution induced by $q$ is stationary,
  let us assume that
  \begin{equation*}\label{eq:niciII}
    \textrm{Prob}(X^0_{\alpha\beta\gamma}(t)=\omega)
    \;=\frac{q(\omega)}{Z^0_n},
    \qquad\textrm{for all }\omega\in\Omega^0_n,
  \end{equation*}
  and try to compute
  $\textrm{Prob}(X^0_{\alpha\beta\gamma}(t+1)=\omega')$.
  For this, recall that $I(t)$ denotes the random wall selected at
  time $t$ and define $J(t+1)$ as follows: if $I(t)$ becomes active so
  that $X^0_{\alpha\beta\gamma}(t+1)=T(X^0_{\alpha\beta\gamma}(t),I(t))$,
  then define $J(t+1)$ by $\bar
  T(X^0_{\alpha\beta\gamma}(t),I(t))=(X^0_{\alpha\beta\gamma}(t+1),J(t+1))$;
  otherwise set $J(t+1)=I(t)$. Then, since $T$ is given as the first
  component of $\bar T$,
  \[
  \textrm{Prob}\big(X^0_{\alpha\beta\gamma}(t+1)=\omega'\big)\;=\;
  \sum_{j=0}^n\textrm{Prob}\big(X^0_{\alpha\beta\gamma}(t+1)=\omega',\;
  J(t+1)=j\big).
  \]
  Now, by definition of the Markov chain $X^0_{\alpha\beta\gamma}$, for
  all $\omega'$ and $j$,
  \begin{eqnarray*}
    \textrm{Prob}\big(X^0_{\alpha\beta\gamma}(t+1)=\omega',\;J(t+1)=j\big)
    &=&
    \lambda(i)\cdot
    \textrm{Prob}\big(X^0_{\alpha\beta\gamma}(t)=\omega,\;I(t)=i\big)\\
    &&\quad\;+\;
    (1-\lambda(j))\cdot\textrm{Prob}\big(X^0_{\alpha\beta\gamma}(t)=\omega',\; 
    I(t)=j\big),
  \end{eqnarray*}
  where $(\omega,i)=\bar T^{-1}(\omega',j)$. Since the random variable
  $I(t)$ is uniform on $\{0,\ldots,n\}$, we get
  \begin{eqnarray*}
    \textrm{Prob}\big(X^0_{\alpha\beta\gamma}(t+1)=\omega',\;J(t+1)=j\big)
    &=&
    \lambda(i)\cdot
    \frac{q(\omega)}{Z^0_n}\frac1{n+1}
    \;+\;
    (1-\lambda(j))\cdot
    \frac{q(\omega')}{Z^0_n}\frac1{n+1}.
  \end{eqnarray*}
  But since $\bar T$ preserves the weights via
  Relation~(\ref{for:preserve}),
  $\lambda(i)q(\omega)=\lambda(j)q(\omega')$ and the terms involving
  $\lambda$ cancel. Finally
  \[
  \textrm{Prob}(X^0_{\alpha\beta\gamma}(t+1)=\omega')\;=\;\sum_{j=0}^n\frac{q(\omega')}{Z^0_n}\frac1{n+1}
  \;=\;\frac{q(\omega')}{Z^0_n},
  \]
  and this completes the proof that the distribution induced by $q$ is
  stationary.
\end{proof}

\section{The 3-TASEP}\label{sec:3tasep}
\begin{figure}[h]
\begin{center}
\includegraphics[width=.5\linewidth]{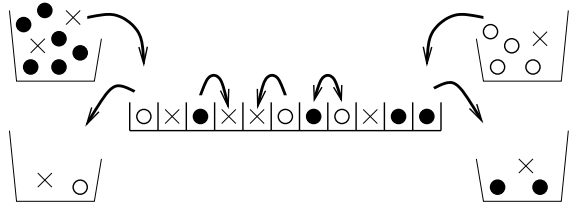}
\end{center}\vspace{-2mm}
\caption{The 3-TASEP.}
\label{fig:rough-extended}
\end{figure}

The combinatorial approach we developed in the previous sections can
be extended to a slightly more general model, the 3-TASEP, which we now
define. The 3-TASEP is similar to the TASEP but each time a black or a
white particle exits, there is a certain probability $\varepsilon$
that the particle that enter in its place is a neutral particle
$\times$. On the one hand, as in the TASEP, black particles always
travel from left to right and white particles always do the opposite.
On the other hand, neutral particles have no preferred direction and
get displaced in opposite direction by black and white particles. An
informal illustration of the 3-TASEP is given by
Figure~\ref{fig:rough-extended}.

\subsection{Definition of the 3-TASEP} 
\begin{figure}[h]
\begin{center}
\includegraphics[width=60mm]{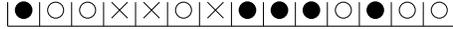}
\end{center}\vspace{-3mm}
\caption{A 3-TASEP configuration  with $n=14$ cells.}
\label{fig:example-config-bwv}
\end{figure}

A \emph{3-TASEP configuration} is a row of $n$ cells, each containing
one particle, which can be of type $\bullet$, $\times$ or $\circ$. An
example of configuration is given by
Figure~\ref{fig:example-config-bwv}.  
The local configuration around a wall $i$ in a configuration $\tau$ is
denoted $\tau[i]$: for $i\in\{1,\ldots,n-1\}$, $\tau[i]$ is the
element of the set
$\{\bullet\act\cross,\bullet\act\circ,\bullet\act\bullet,\cross\act\bullet,\cross\act\cross,\cross\act\circ,\circ\act\bullet,\circ\act\cross,\circ\act\circ\}$
that describes the two cells surrounding wall $i$, for $i=0$,
$\tau[0]\in\{\act\bullet,\act\cross,\act\circ\}$, and for $i=n$,
$\tau[n]\in\{\bullet\act,\cross\act,\circ\act\}$.

The 3-TASEP is a Markov chain $S_{\alpha\beta\gamma\varepsilon}$
defined on the set of 3-TASEP configurations in terms of four
parameters $\alpha$, $\beta$ and $\gamma$ in $]0,1]$ and $\varepsilon$
in $[0,1]$. From time $t$ to $t+1$, the chain evolves from the
configuration $\tau=S_{\alpha\beta\gamma\varepsilon}(t)$ to a
configuration $\tau'=S_{\alpha\beta\gamma\varepsilon}(t+1)$ as follows:
\begin{itemize}
\item A wall $i$ is chosen uniformly at random among the $n+1$ walls.
\item Depending on the local configuration $\tau[i]$ around wall $i$, a
transition may be triggered:
\begin{itemize}
\item unstable local configurations in the middle ($i\in\{1,\ldots,n-1\}$):
\begin{itemize}
\item[$a_1$.] Case $\bullet\act\circ$, a transition
$\bullet\act\circ\to\circ|\bullet$ occurs with probability
$\lambda(\bullet\act\circ):=\alpha$.
\item[$a_2$.] Case $\cross\act\circ$, a transition
$\cross\act\circ\to\circ|\cross$ occurs with probability
$\lambda(\cross\act\circ):=\beta$.
\item[$a_3$.] Case $\bullet\act\cross$, a transition
$\bullet\act\cross\to\cross|\bullet$ occurs with probability
$\lambda(\bullet\act\cross):=\gamma$.
\end{itemize}
\item unstable local configurations on the left border ($i=0$):
\begin{itemize}
\item[$b_1$.] Case $\act\circ$, the  particle exits with total
probability $\lambda(\act\circ):=\beta$, in 2 possible ways:
\begin{itemize}
\item[$b_1'$.] a transition $\act\circ\to|\bullet$ occurs
with probability
$(1-\varepsilon)\,\beta,$
\item[$b_1''$.] or a transition $\act\circ\to|\cross$ with
probability $\varepsilon\beta$
(neutralization),
\end{itemize}
\item[$b_2$.] Case $\act\cross$, a transition $\act\cross\to|\bullet$
occurs with probability
$\lambda(\act\cross):=
(1-\varepsilon)\beta\,\gamma/\alpha$.
\end{itemize}
\item unstable local configurations on the right border ($i=n$):
\begin{itemize}
\item[$c_1$.] Case $\bullet\act$, the particle exits with total
probability $\lambda(\bullet\act):=\gamma$, in 2 possible ways:
\begin{itemize}
\item[$c_1'$.] a transition $\bullet\act\to\circ|$ occurs with
probability $(1-\varepsilon)\,\gamma$,
\item[$c_1''$.] or a transition $\bullet\act\to\cross|$ with
probability $\varepsilon\gamma$
(neutralization),
\end{itemize}
\item[$c_2$.] Case $\cross\act$, a transition $\cross\act\to\circ|$
occurs with probability
$\lambda(\cross\act):=(1-\varepsilon)\gamma\,\beta/\alpha$.
\end{itemize}
\item stable local configurations:
\begin{itemize}
\item[$d$.] Cases $\bullet\act\bullet$, $\cross\act\cross$,
  $\circ\act\circ$, $\circ\act\cross$, $\circ\act\bullet$,
  $\cross\act\bullet$, $\act\bullet$ and $\circ\act$, no 
  transition occur:
  $\lambda(\bullet\act\bullet)=\lambda(\cross\act\cross)
  =\lambda(\circ\act\circ)=\lambda(\circ\act\cross)
  =\lambda(\circ\act\bullet)=\lambda(\cross\act\bullet)
  =\lambda(\act\bullet)=\lambda(\circ\act):=0$.
\end{itemize}
\end{itemize}
\item If a transition occurs, the new configuration $\tau'$ is
obtained from $\tau$ by applying the transition to the local
configuration around the chosen wall. Otherwise, $\tau'=\tau$.
\end{itemize}
In order to explain the role of the parameters $\alpha$, $\beta$,
$\gamma$ and $\varepsilon$ a few remarks are useful:
\begin{itemize}
\item The equality $\lambda(\cross\act\circ)=\lambda(\act\circ)$
  translates the idea that a white particle feels the same attraction
  to the left in front of a neutral particle as it feels for exiting
  at the left border.  A similar interpretation holds for
  $\lambda(\bullet\act\cross)=\lambda(\bullet\act)$ and black particles.
\item The equality $\lambda(\cross\act)/\lambda(\cross\act\circ)
  =(1-\varepsilon)\lambda(\bullet\act)/\lambda(\bullet\act\circ)$, or
  says that the ratio between entry and movement rates for white
  particles is the same in presence of black or neutral particles.  A
  similar interpretation holds for black particles.
\item The fact that the same quantity $\varepsilon$ controls the
  probability that a $\times$ particles enters instead of a black
  particle or instead of a white particle may be thought of as a
  curious restriction: it is dictated by technical considerations in
  the proof, and at the present state we do not know whether it can be
  easily circumvented or not.
\end{itemize}

The TASEP with parameter $\alpha$, $\beta$ and $\gamma$ is recovered
by taking $\varepsilon=0$. Indeed, in this case, after the initial
neutral particles have exit the system, no new neutral particles are
created and the rules are exactly those of the TASEP as presented in
Section~\ref{JumpingGilles}.

It will be useful to reformulate again the transition of the 3-TASEP
in terms of applications from the set of configurations with a chosen
wall into the set of configurations. Since there are two possible
transitions in the cases $\act\circ$ and $\bullet\act$ we introduce the
following two applications:
 \begin{itemize}
 \item The application $\vartheta_1:(\tau,i)\to\tau'$ performing at
   wall $i$ the transitions prescribed by cases $a_1$, $a_2$ and
   $a_3$, $b'_1$ and $b_2$, $c'_1$ and $c_2$.
 \item The application $\vartheta_2:(\tau,i)\to\tau'$ performing at
   wall $i$ the transitions prescribed by cases $a_1$, $a_2$ and
   $a_3$, $b''_1$ and $b_2$, $c''_1$ and $c_2$.
 \end{itemize}
 Then the transitions of the chain $S_{\alpha\beta\gamma\varepsilon}$
 can be described as follows: choose $i=I(t)$ uniformly at random in
 $\{0,\ldots,n\}$ and set
\[
S_{\alpha\beta\gamma\varepsilon}(t+1)=\left\{
\begin{array}{llc}
\vartheta_1(\tau,i)
& \textrm{ with probability} &(1-\varepsilon) \lambda(\tau[i]),\\
\vartheta_2(\tau,i)
& \textrm{ with probability} &\varepsilon \lambda(\tau[i]),\\
\tau & 
\textrm{ otherwise,}
\end{array}\right.
\]
where $\tau=S_{\alpha\beta\gamma\varepsilon}(t)$, and $\tau[i]$
denotes the local configuration around wall $i$ in $\tau$.

\subsection{Complete configurations for the 3-TASEP}
The complete configurations for the 3-TASEP are concatenations of
complete configurations for the TASEP separated by $|\xx|$-columns:
more explicitly, each complete configuration $\omega$ with $\ell$
$\times$-particles in the first row can be uniquely written
$\omega_0|\xx|\omega_1\cdots|\xx|\omega_\ell$ where each $\omega_i$
belongs to $\Omega_{n_i}$ for some $n_i\geq0$.  In other terms, these
complete configurations are pairs of rows of cells containing
particles such that the $\times$-particles always form $|\xx|$-columns
and such that between two $|\xx|$-columns the balance and positivity
conditions are satisfied. Let $\Omega_n$ denote the set of
complete configurations of length $n$.

\begin{figure}[h]
  \begin{center}
    \includegraphics[width=60mm]{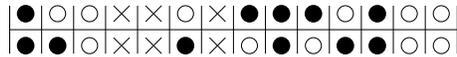}
  \end{center}
  \vspace{-3mm}
  \caption{A complete configuration  with $n=14$.}
  \label{fig:example-complete-config-x}
\end{figure}
An example of a complete configuration is given in
Figure~\ref{fig:example-complete-config-x}: from left to right the
subconfigurations have successively length 3, 0, 1, and~7.  The local
configuration around wall $i$ in a complete configuration $\omega$,
describing the one or two columns surrounding wall $i$, is denoted
$\omega[i]$. The following enumerative Lemmas are proved in
Section~\ref{sec:counting}.

\begin{lemma} \label{lem:count-all} 
  The cardinality of
  $\Omega_n$ is $\frac{1}{2}{2n+2\choose n+1}$.
\end{lemma}
\begin{lemma} \label{lem:count-all-refined}
  For any $k+\ell+m=n$, the cardinality of the set $\Omega_{k,m}^\ell$
  of complete configurations of $\Omega_n$ with $\ell$
  $|\xx|$-columns, and $k$ black and $m$ white particles on the top
  row is $\frac{\ell+1}{n+1}{n+1\choose k}{n+1\choose m}$.
\end{lemma}
\begin{lemma} \label{lemma:olaf-le-retour}
  For any $\ell+p=n$, the cardinality of the set $\Omega_{n}^\ell$ of
  complete configurations of $\Omega_n$ with $\ell$ $|\xx|$-columns is
  $\frac{\ell+1}{n+1}{2n+2\choose n-\ell}$.
\end{lemma}

\subsection{The complete chain $X_{\alpha\beta\gamma\varepsilon}$} 
We shall directly define the chain $X_{\alpha\beta\gamma\varepsilon}$
in terms of two bijections $\bar T_1$ and $\bar T_2$ from
$\Omega_n\times\{0,\ldots,n\}$ to itself.
To do this, we partition the set $\Omega_n\times\{0,\ldots,n\}$ into
classes, and we first describe for each class $A$ the image
$(\omega',j)$ of a pair $(\omega,i)\in A$ by $\bar T_1$. The bijection
$\bar T_2$ is then described as a simple variation on $\bar T_1$.

As in Section~\ref{sec:informal}, given a complete configuration
$\omega$ with top row $\tau$ and a wall $i$, we distinguish the
following walls: if the local configuration $\tau[i]$ is
$\bullet\act\circ$, $\cross\act\circ$, $\bullet\act$ or $\cross\act$,
then let $j_1<i$ be the leftmost wall such that there are only white
particles in the top row between walls $j_1$ and $i-1$; if the local
configuration $\tau[i]$ is $\bullet\act\circ$, $\bullet\act\cross$,
$\act\circ$ or $\act\cross$, then let $j_2>i$ be the rightmost wall
such that there are only black particles in the top row between walls
$i+1$ and $j_2$.

\begin{figure}
  \begin{center}
    \includegraphics[width=\linewidth]{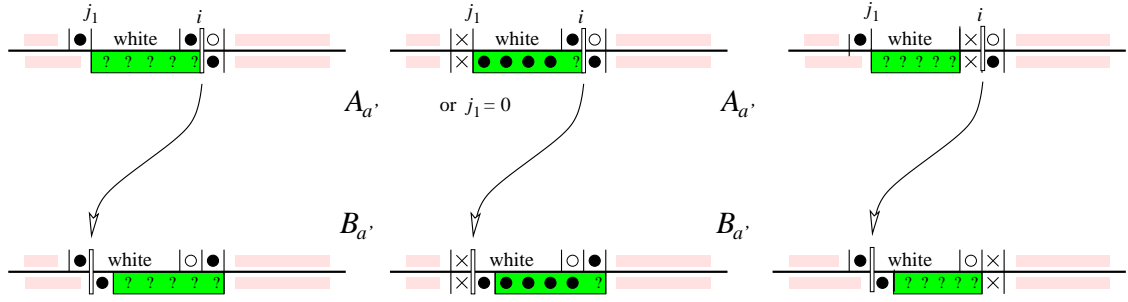}
  \end{center}
  \caption{Cases $\bullet\act\circ$ and $\times\act\circ$ for the bijections $\bar T_1$ and $\bar T_2$.}
  \label{move-right-sweep-x}
\end{figure}

The action of $\bar T_1$ is described separately for the different
cases of local configuration $\omega[i]$:
\begin{itemize} 
\item unstable local configurations in the middle
  ($i\in\{1,\ldots,n-1\}$):
  \begin{itemize}
  \item[$A_{a'}$.] Cases ${\bullet\atop?}\act\wb$ and
    $\xx\act{\circ\atop\bullet}$. Then $j=j_1$, and $\omega'$ is obtained
    by moving the $|\wb|$-column from the right-hand side of wall $i$ to the
    right-hand side of wall $j$ (Figure~\ref{move-right-sweep-x}).
    
    The image $B_{a'}$ of this class consists of pairs $(\omega',j)$
    such that: the  wall $j$ is the left border ($j=0$) or there is a
    black or a $\times$ particle on its left-hand side, there is a
    $|\wb|$-column on its right-hand side, and the sequence of white
    particles on the right-hand side of wall $j$ in the top row is
    followed by a black or a $\times$ particle.

  \item[$A_{a''}$.] Cases ${\bullet\atop?}\act\ww$ or
    ${\bullet\atop\circ}\act\xx$. Then $j=j_2$, and $\omega'$ is
    obtained by removing the two particles that form the
    ${\bullet\atop}|{\atop\circ}$-diagonal or the $|\bw|$-column at
    wall $i$ and replacing them at wall $j$ so that they form a
    ${\bullet\atop}|{\atop\circ}$-diagonal if there is a white
    particle on the right-hand side of wall $j$ in the top row, or a
    $|\bw|$-column otherwise (Figure~\ref{move-left-sweep-x}).
    
    The image $B_{a''}$ of this class consists of pairs
    $(\omega',j)$ with a $|\ww|$-column, an $|\xx|$-column, or the
    border on the right-hand side of wall $j$ and such that there is a
    non-empty sequence of black particles on the left-hand side of
    wall $j$ in the top row, followed by a white or an $\times$
    particle.
  \end{itemize}    

\begin{figure}
  \begin{center}
    \includegraphics[width=\linewidth]{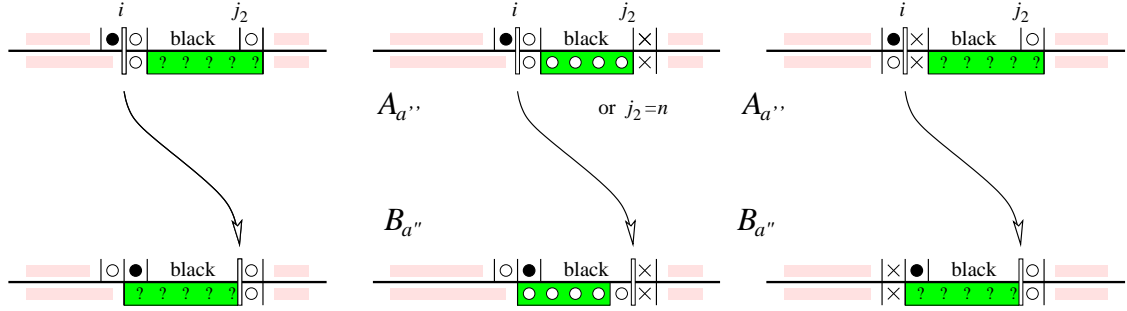}
  \end{center}
  \caption{Cases $\bullet\act\circ$ and $\bullet\act\times$ 
    with black sweep for the bijections
    $\bar T_1$ and $\bar T_2$.}
  \label{move-left-sweep-x}
\end{figure}

\item unstable local configurations on the left border $(i=0$):
  \begin{itemize}
  \item[$A_{b'}$.] Case $\act\wb$. Then $j=j_2$, and $\omega'$ is
    obtained by removing the two particles that form the
    $|\wb|$-column on the left border and replacing them at wall $j$
    so that they form a ${\bullet\atop}|{\atop\circ}$-diagonal if
    there is a white particle on the right-hand side of wall $j$ in
    the top row, or a $|\bw|$-column otherwise.

    The image $B_{b'}$ of the class $A_{b'}$ by $\bar T_1$ consists of
    pairs $(\omega',j)$ with a $|\ww|$-column, an $|\xx|$-column, or
    the border on the right-hand side of wall $j$ of $\omega'$ and such
    that there is a non-empty sequence of black particles on the
    left-hand side of wall $j$ in the top row, ending at the left
    border.
  \item[$A_{b''}$.] Case $\act\xx$. Then $\omega'=\omega$ and $j=0$.  
    The image $B_{b''}$ of the class $A_{b''}$ by $\bar T_1$ consists of
    pairs $(\omega',0)$ with a $|\xx|$-column on the left border.
  \end{itemize}
\item unstable local configurations on the right border ($i=n$):
  \begin{itemize}
  \item[$A_{c'}$] Case $\bw\act$. Then $j=j_1$, and $\omega'$ is
    obtained by removing the rightmost $|\bw|$-column and forming a
    $|\wb|$-column on the right-hand side of wall $j$.
    
    The image $B_{c'}$ of the class $A_{c'}$ by $\bar T_1$ consists of
    pairs $(\omega',j)$ such that: the wall $j$ is the left border
    ($j=0$) or it has a black or a $\times$ particle on its left-hand
    side, there is a $|\wb|$-column on its right-hand side, and the
    sequence of white particles on the right-hand side of wall $j$ in
    the top row ends at the right border.
  \item[$A_{c''}$.] Case $\xx\act$. Then $\omega'=\omega$ and $j=n$. The
    image $B_{c''}$ of the class $A_{c''}$ by $\bar T_1$ consists of pairs
    $(\omega',n)$ with a $|\xx|$-column on the right border.
  \end{itemize}
\end{itemize}

\begin{figure}  
  \begin{center}
    \includegraphics[width=.95\linewidth]{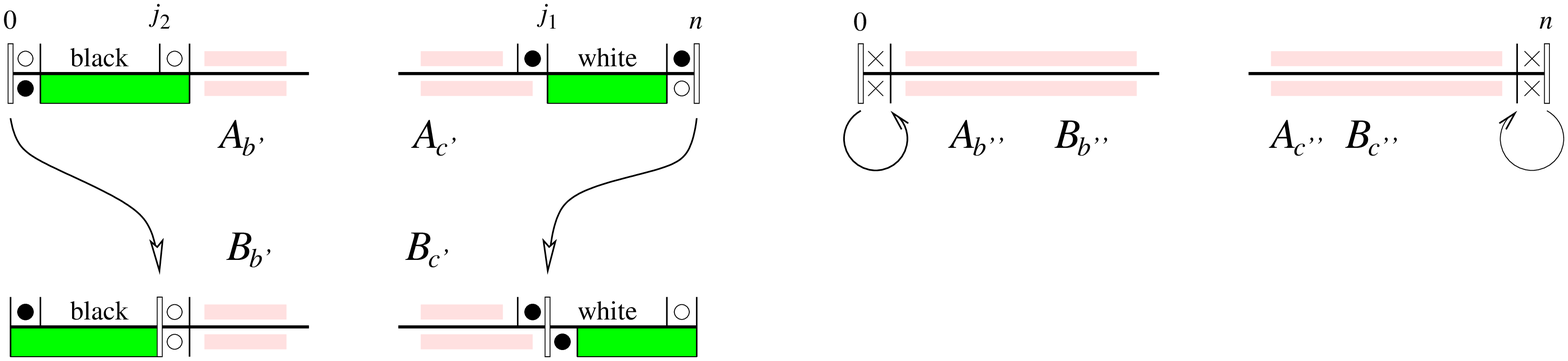}
  \end{center}
  \caption{The application $\bar T_1$ on the borders.}
  \label{fig:enter-T1}
  
  \smallskip
  
  \begin{center}
    \includegraphics[width=.95\linewidth]{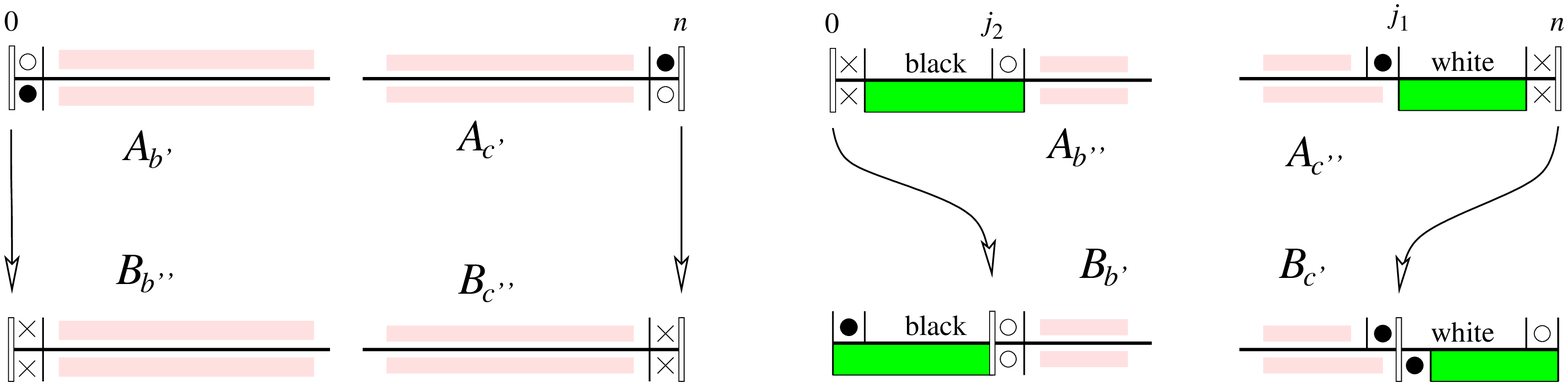}
  \end{center}
  \caption{The application $\bar T_2$ on the borders.}
  \label{fig:enter-T2}
\end{figure}

Finally let $A_d$ denote the set of pairs $(\omega,i)$ such that the
local configuration around wall $i$ of $\omega$ is stable. The
mapping $\bar T_1$ has no effect on these pairs, and $B_d=A_d$.

The application is invertible in each case and the sets
$\{B_{a'},B_{a''},B_{b'},B_{b''},B_{c'},B_{c''},B_d\}$ form a
partition of $\Omega_n\times\{0,\ldots,n\}$. Hence $\bar T_1$
is a bijection from $\Omega_n\times\{0,\ldots,n\}$ onto itself.

The application $\bar T_2$ differs from $\bar T_1$ only at the
borders. Consider the involution $Y$ on $\Omega_n\times\{0,\ldots,n\}$
that acts only on a pair $(\omega,i)$ by changing, if $i=0$ or $i=n$,
the local configuration $\omega[i]$ according to the following rules:
\[\textstyle
\act\wb\;\leftrightarrow\;\act\xx, \qquad\textrm{ and }\quad
\bw\act\;\leftrightarrow\;\xx\act.
\]
Then the image of $(\omega,i)$ by $\bar T_2$ is defined to be the
image of $Y(\omega,i)$ by $\bar T_1$. In particular $\bar T_2$, being
the composition $\bar T_1\circ Y$ of two bijections is itself a
bijection. The action of $\bar T_2$ on the borders is illustrated by
Figure~\ref{fig:enter-T2}.

Let now $T_1$ and $T_2$ denote respectively the first component of
$\bar T_1$ and $\bar T_2$. Then the Markov chain
$X_{\alpha\beta\gamma\varepsilon}$ is defined in terms of the
$T_1$ and $T_2$ exactly as $S_{\alpha\beta\gamma\varepsilon}$
is defined in terms of the $\vartheta_1$ and
$\vartheta_2$: choose $i=I(t)$ uniformly at random in $\{0,\ldots,n\}$
and set
\[
X_{\alpha\beta\gamma\varepsilon}(t+1)=\left\{
\begin{array}{llc}
T_1(\omega,i)
& \textrm{ with probability } & (1-\varepsilon) \lambda(\tau[i]),\\
T_2(\omega,i)
& \textrm{ with probability }& \varepsilon \lambda(\tau[i]),\\
\omega & 
\textrm{ otherwise,}
\end{array}\right.
\]
where $\omega=X_{\alpha\beta\gamma\varepsilon}(t)$, and $\tau=\textrm{top}(\omega)$.

\subsection{The stationary distribution 
  of $X_{\alpha\beta\gamma\varepsilon}$ and
  $S_{\alpha\beta\gamma\varepsilon}$}

The parameter $n_y$ and $n_z$ of Section~\ref{sec:arbitrary} are
extended in a straightforward way to complete configurations of
$\Omega_n$ by putting labels independently in each subconfiguration
delimited by $|\xx|$-columns or borders. Then for
$\omega\in\Omega_n$, set
\[
q(\omega)\;=\;
\beta^n\gamma^n(1-\varepsilon)^n
\left(\frac{\alpha}{\beta}\right)^{n_y(\omega)}
\left(\frac{\alpha}{\gamma}\right)^{n_z(\omega)}
\left(\frac{\alpha^2\varepsilon}
{\beta\gamma(1-\varepsilon)}\right)^{\ell(\omega)}
.
\]
where $\ell(\omega)$ denotes the number of $|\xx|$-columns in $\omega$.

Then Theorem~\ref{thm:generic} extends verbatim:
\begin{theorem}\label{thm:generic-3tasep}
  The Markov chain $X_{\alpha\beta\gamma\varepsilon}$ has the
  following unique stationary distribution:
  \[
  \textrm{Prob}(X_{\alpha\beta\gamma\varepsilon}(t)=\omega)
  \;\mathop{\longrightarrow}_{t\rightarrow\infty}\;
  \frac{q(\omega)}{Z_n}, \qquad
  \textrm{where } Z_n={\sum_{\omega' \in \Omega_n}q(\omega')},
  \]
  where $q(\omega)$ is the previously defined weight on the complete
  configurations of the 3-TASEP.
\end{theorem}

Again this theorem immediately yields a combinatorial interpretation
of the stationary distribution of the chain
$S_{\alpha\beta\gamma\varepsilon}$, via the relation
$S_{\alpha\beta\gamma\varepsilon}=
\mathrm{top}(X_{\alpha\beta\gamma\varepsilon})$.  In particular
in the case $\alpha=\beta=\gamma=1$, $\varepsilon=1/2$, we
obtain the following corollary on $S=S_{111\frac12}$ and $X=X_{111\frac12}$:
\begin{corollary}
  Let $\textrm{top}(\omega)$ denote the top row of a complete
  configuration $\omega$. Then for any initial distributions
  $S(0)$ and $X(0)$ with
  $\textrm{top}(X(0))=S(0)$, and any basic
  configuration~$\tau$,
 \[
 \textrm{Prob}(S(t) = \tau)\;=\;
 \textrm{Prob}(\textrm{top}(X(t)) =\tau)
 \;\mathop{\longrightarrow}_{t\rightarrow\infty}\;
 \frac{\big|\{\omega\in\Omega_n\mid
 \textrm{top}(\omega)=\tau\}|}{|\Omega_{n}|}.
 \]
In particular, for any $k+\ell+m=n$, we obtain combinatorially the
formula:
\[
\textrm{Prob}(S(t) \textrm{ contains } k \textrm{ black and } m
 \textrm{ white particles})
 \;\mathop{\longrightarrow}_{t\rightarrow\infty}\;
 \frac{|\Omega^\ell_{k,m}|}{|\Omega_n|} \;=\;
 \frac{\frac{\ell+1}{n+1}{n+1\choose k}{n+1\choose m}}{\frac{1}{2}
 {2n+2\choose n+1}}.
 \]
\end{corollary}

Theorem~\ref{thm:generic-3tasep} is an easy consequence of the fact
that the two bijections preserve weights in the sense of the following
lemma. Recall that $\lambda$ describes the transition
probabilities for each possible local configuration. 
\begin{lemma}\label{lem:preserve-x}
  The applications $\bar T_1$ and $\bar T_2$ transport together the
  weight $\lambda$ in the following sense: for all
  $(\omega',j)\in\Omega\times\{0,\ldots,n\}$,
\[
(1-\varepsilon)\lambda(\omega_1[i_1])\,q(\omega_1)
+\varepsilon\lambda(\omega_2[i_2])\,q(\omega_2)
\;=\;\lambda(\omega'[j])\,q(\omega'),
\]
where $(\omega_1,i_1)=\bar T_1^{-1}(\omega',j)$ and $(\omega_2,i_2)=\bar
T_2^{-1}(\omega',j)$.
\end{lemma}
\begin{proof}
  This lemma is easily verify by a case by case analysis similar to
  that of Lemma~\ref{lem:preserve}.
\end{proof}

\begin{proof}[Proof of Theorem~\ref{thm:generic-3tasep}]
This proof exactly mimics the proof of Theorem~\ref{thm:generic},
using Lemma~\ref{lem:preserve-x} instead of Lemma~\ref{lem:preserve}.
\end{proof}

\section{Periodic boundary conditions}\label{sec:periodic}
\begin{figure}
  \begin{center}
    \includegraphics[width=.5\linewidth]{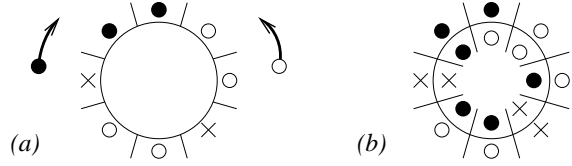}
  \end{center}
  \caption{A basic and a complete configuration of the 3-TASEP with
  periodic boundary conditions.}
  \label{fig:informal-circle}
\end{figure}

A standard alternative to our definition of the TASEP is to consider
periodic boundary conditions: the leftmost cell is considered on the
right-hand side of the rightmost cell, or equivalently, the
configurations are arranged on a circle (see
Figure~\ref{fig:informal-circle} a., the circle is rigid, not subject
to rotation).

Since there are no border walls in these configurations, the Markov
chain $\smash{\widehat S_{\alpha\beta\gamma}}$ is defined using only
Cases $a_1$, $a_2$, $a_3$ of the transition of the 3-TASEP.  Observer
that the numbers $k$, $\ell$ and $m$ of black, $\cross$ and white
particles do not change during the evolution. The case without
$\times$ particle is easily seen to have a uniform stationary
distribution, so we concentrate on the case with at least one
$\times$ particle.

Our approach is easily adapted to deal with this case.  Let
$\widehat\Omega_n$ be a new set of complete configurations that are
made of two rows of cells arranged on a circle and that are such that
the subconfigurations between two $|\xx|$-columns, when read in
clockwise direction, satisfy the balance and positivity constraints.
More precisely we are interested in the subset
$\smash{\widehat\Omega^\ell_{k,m}}$ of configurations of
$\smash{\widehat\Omega_n}$ that have $\ell$ $|\xx|$-columns, $k$ black
and $m$ white particles in the top row. The following lemma is proved
in Section~\ref{sec:counting}.
\begin{lemma}\label{lemma:nici}
The cardinality of $\widehat \Omega^\ell_{k,m}$ is
${{n \choose k}{n \choose m}}$.  
\end{lemma}
Cases~$A_{a'}$ and $A_{a''}$ of the definition of $\bar T_1$ allow to
define a bijection $\widehat T$ from
$\smash{\widehat\Omega^\ell_{k,m}\times\{0,\ldots,n-1\}}$ to itself
and an associated Markov chain $\widehat X_{\alpha\beta\gamma}$ such
that $\widehat S_{\alpha\beta\gamma}\equiv\mathrm{top}(\widehat
  X_{\alpha\beta\gamma})$.  The same argument as in
  Section~\ref{sec:maximal} for the chain $X^0$ then immediately
  yields the fact that $\widehat X=\widehat X_{111}$ has a uniform
  stationary distribution. In particular:
\[
\textrm{Prob}(\widehat X(t)=\omega)
\;\mathop{\longrightarrow}_{t\rightarrow\infty}\;
\frac1{|\widehat\Omega^\ell_{k,m}|}\;=\;\frac1{{n \choose k}{n \choose
    m}}.
\]

Furthermore, the statistics $n_y$ and $n_z$ are immediately extended
to configurations of $\widehat\Omega_n$ by putting label independently
on every subconfiguration between $|\xx|$-columns.
Lemma~\ref{lem:preserve-x} adapts in a straightforward way (with
$\varepsilon=0$), and allows to express the stationary distribution in the
general case:
\begin{theorem}
  The Markov chain $\widehat X_{\alpha\beta\gamma}$ has the following
  unique stationary distribution:
\[
\textrm{Prob}(\widehat X_{\alpha\beta\gamma}(t)=\omega)
\mathop{\longrightarrow}_{t\to\infty}\frac{q(\omega)}{\widehat Z_n},
\qquad\textrm{ where }\widehat Z_n=\sum_{\omega'\in\widehat
  \Omega_n}q(\omega'),
\]
where
$q(\omega)=\beta^n\gamma^n(\alpha/\beta)^{n_y(\omega)}(\alpha/\gamma)^{n_z(\omega)}$.
\end{theorem}

Finally the stationary distribution of $S_{\alpha\beta\gamma}$ is
recovered from the relation $\widehat S_{\alpha\beta\gamma}\equiv
\widehat X_{\alpha\beta\gamma}$. In particular,
for $\widehat S=\widehat S_{111}$, 
\[
\textrm{Prob}(\widehat S(t)=\tau)
\;\mathop{\longrightarrow}_{t\rightarrow\infty}\;
\frac{|\{\omega\in\widehat\Omega^\ell_{k,m}
  \mid\textrm{top}(\omega)=\tau\}|} {|\widehat\Omega^\ell_{k,m}|}.
\] 
For instance, a configuration $\tau$ of the form
\[|\cross|
\underbrace{\circ\cdots\circ}_{m_1}|
\underbrace{\bullet\cdots\bullet}_{k_1}|
\cross|
\underbrace{\circ\cdots\circ}_{m_2}|
\underbrace{\bullet\cdots\bullet}_{k_2}|
\cross|\cdots|
\cross|
\underbrace{\circ\cdots\circ}_{m_\ell}|
\underbrace{\bullet\cdots\bullet}_{k_\ell}|
\]
 for some $k_1+\ldots+k_\ell=k$,
$m_1+\ldots+m_\ell=m$ corresponds to only one complete configuration
\[\textstyle
|\xx|
\underbrace{\textstyle\wb\cdots\wb}_{m_1}|
\underbrace{\textstyle\bw\cdots\bw}_{k_1}|
\xx|
\underbrace{\textstyle\wb\cdots\wb}_{m_2}|
\underbrace{\textstyle\bw\cdots\bw}_{k_2}|
\xx|\cdots|
\xx|
\underbrace{\textstyle\wb\cdots\wb}_{m_\ell}|
\underbrace{\textstyle\bw\cdots\bw}_{k_\ell}|
\]
(because of the positivity constraints on blocks between
$|\xx|$-columns), and thus has probability 
$1/{n\choose k}{n\choose m}$
in the stationary distribution of $\widehat S$.

\section{Irreducibility}\label{sec:paths}

In this section we verify that the Markov chains $X^0$, $\widehat X$
and $X$ are irreducible, \emph{i.e.} that there is a positive
probability to go from any configuration $\omega$ to any other one
$\omega'$.  In other terms we need to prove that the transition graphs
of these three chains are connected. The proof is based on an
observation about iterating the bijections $\bar T$, or $\bar T_1$ or
$\bar T_2$, and on induction on $n$.

To every pair $(\omega,i)$ of $\Omega_n\times\{0,\ldots,n\}$ we
associate a reduced configuration $\omega^i$ in $\Omega_{n-1}$,
obtained from $\omega$ by deleting two particles around the wall $i$
in a way that depends on the local configuration:
\begin{itemize}
\item Cases ${\bullet\atop?}\act\wb$, $\xx\act\wb$ and $\act\wb$.  The
  reduced configuration $\omega^i$ is obtained by removing the
  $|\wb|$-column on the right-hand side of wall $i$.
\item Case ${\bullet\atop?}\act\ww$. The reduced configuration $\omega^i$ is obtained by removing the two particles forming the ${\bullet\atop}|{\atop\circ}$-diagonal around wall $i$.
\item Cases $\bw\act\xx$ and $\bw\act$. The reduced configuration
  $\omega^i$ is obtained by removing the $|\bw|$-column on the
  left-hand side of wall $i$.
\item Cases $\act\xx$ and $\xx\act$. The reduced configuration is
  obtained by removing the $|\xx|$-column on the border.
\end{itemize}

\begin{lemma}\label{lem:cycle}
  Let $\tilde\omega$ be a configuration of $\Omega_{n-1}$. Let
  $S(\tilde\omega)$ be the set of pairs $(\omega,i)$ of
  $\Omega_n\cross\{0,\ldots,n\}$ having $\tilde\omega$ as reduced
  configuration, \emph{i.e.} such that $\omega^i=\tilde\omega$. In
  particular let $\omega_0$ be the configuration $|\xx|\tilde\omega$ and
  $\omega_n$ be the configuration $\tilde\omega|\xx|$, and define
  $S^0(\tilde\omega)=S(\tilde\omega)\setminus\{(\omega_0,0),(\omega_n,n)\}$.
  Then:
  \begin{itemize} 
  \item The set $S^0(\tilde\omega)$ is a cyclic orbit of $\bar T_1$:
    given $(\omega,i)\in S^0(\tilde\omega)$, all other elements of
    $S^0(\tilde\omega)$ can be reached by successive applications of
    $\bar T_1$.
  \item The set $S(\tilde\omega)$ is a cyclic orbit of $\bar T_2$.
  \item If $\tilde\omega\in\Omega^0_{n-1}$ then
    $S^0(\tilde\omega)\subset\Omega^0_n$ and $S^0(\tilde\omega)$ is a
    cyclic orbit of $\bar T$.
  \end{itemize}
\end{lemma}
\begin{proof}
  As can be checked on Figures~\ref{move-right-sweep}
  and~\ref{move-right-sweep-x}, starting from a pair $(\omega,i)$ of
  the corresponding classes and iterating $\bar T_1$, $\bar T_2$ or
  $\bar T$, the selected wall moves to the left with the pair of
  black and white particles, and successively stops on the right-hand
  side of every black or $\cross$ particle of the top row, until it
  reaches the left border.  Similarly, as can be checked on
  Figures~\ref{move-left-sweep} and~\ref{move-left-sweep-x}, iterating
  $\bar T_1$, $\bar T_2$ or $\bar T$ from a pair $(\omega,i)$ of the
  corresponding classes, the selected wall moves to the right with the
  pair of black and white particles, stopping on the left-hand side of
  every white and $\cross$ particle of the top row, until it reaches
  the right border.
  
  As shown by Figures~\ref{fig:enter-T1}--\ref{fig:enter-T2}, the
  application $\bar T_2$, and the applications $\bar T_1$ or $\bar T$
  behave differently when the border is reached: $\bar T_2$ visits the
  configurations $\omega_0$ or $\omega_n$ while $\bar T_1$ or $\bar
  T$ skips them and immediately restart moving in the opposite
  direction.
  
  Starting from an element $(\omega,i)$ all other elements of
  $S(\tilde\omega)$ (respectively $S^0(\tilde\omega)$) are
  thus visited in a cycle by successive applications of $\bar T_2$
  (respectively $\bar T_1$ or $\bar T$).
\end{proof}
Lemma~\ref{lem:cycle} provides us with cycles in the transition graph
on $\Omega_n$, and each cycle is associated to a reduced configuration
of $\Omega_{n-1}$.  The next lemma transports transitions from
$\Omega_{n-1}$ to $\Omega_n$.
\begin{lemma}\label{lem:releve}
  Let $(\tilde\omega',j)=\bar T_1(\tilde\omega,i)$  be a transition
  between two configurations of $\Omega_{n-1}$. Then there exist
  $k,i_+,j_+$ and $\omega,\omega'$ such that $(\omega,k)\in
  S(\tilde\omega)$, $(\omega',k)\in S(\tilde\omega)$, and $(\omega',j_+)=\bar
  T_1(\omega,i_+)$. The same holds for $\bar T_2$.
\end{lemma}
\begin{proof} For $\bar T_1$ observe that in each case of 
  Figure~\ref{move-right-sweep-x}, and on the second leftmost case of
  Figure~\ref{fig:enter-T1}, a $|\wb|$-column can be inserted on the
  left border without interfering with the action of $\bar T_1$: take
  $\omega=|\wb|\tilde\omega$, $\omega'=|\wb|\tilde\omega$, $k=0$,
  $i_+=i+1$, $j_+=j+1$.  Similarly in each case of
  Figure~\ref{move-left-sweep-x}, and on the leftmost case of
  Figure~\ref{fig:enter-T1}, a $|\bw|$-column can be inserted on the
  right border without interfering with the action of $\bar T_1$: take
  $\omega=\tilde\omega|\bw|$, $\omega'=\tilde\omega|\bw|$, $k=n$,
  $i_+=i$, $j_+=j$.
  
  For $\bar T_2$ observe that in each case of
  Figures~\ref{move-left-sweep-x}--\ref{fig:enter-T2}, an
  $|{\cross\atop\cross}|$-column can be inserted, either on the left
  or on the right border, without interfering with the action of $\bar
  T_2$.
\end{proof}
Lemma~\ref{lem:releve} gives a transition between an element of the
cycle associated to $\tilde\omega$ and an element of the cycle associated to
$\tilde\omega'$. Taking the connectivity of the transition graph on
$\Omega_{n-1}$ as induction hypothesis, we conclude that all cycles of
Lemma~\ref{lem:cycle} belong to the same connected component of the
transition graph defined by $\bar T_2$ on $\Omega_n$. Since every
element of $\Omega_n$ belong to a cycle, this concludes the proof of
the irreducibility of $X$. 

As opposed to this the transition graph defined by $\bar T_1$ is seen
to connect only configurations with the same number of
$|\xx|$-columns. In particular the chain $X_{\alpha\beta\gamma0}$ with
$\varepsilon=0$ is not irreducible, but, the transition graph
defined by $\bar T$ (or $\bar T_1$) on $\Omega^0_n$ is connected and
the chain $X^0_{\alpha\beta\gamma}$ is irreducible. 

Finally the chain $\widehat T$ is seen to be irreducible in a similar
manner as soon as there is at least one $|\xx|$-column.

\section{The number of complete configurations and the cycle
  lemma}\label{sec:counting}

\smallskip
\noindent\textbf{Lemma~\ref{lem:count-all}}. {\em
  The cardinality of $\Omega_n$ is
  $\frac{1}{2} {2n+2\choose n+1}$.}\par\smallskip

\begin{figure}
\begin{center}
\includegraphics[width=120mm]{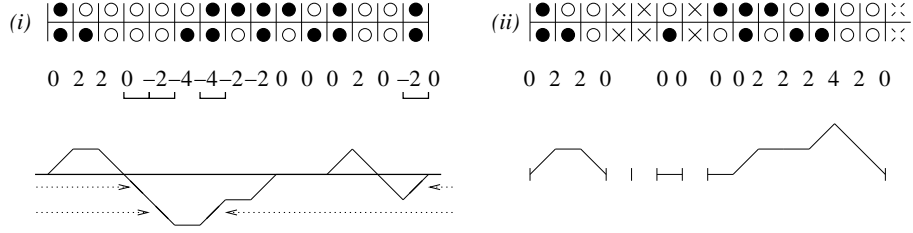}
\end{center}
\caption{From (i) an element of $\overline\Gamma_{n+1}$, to (ii) one
  of $\Omega_n$. The $(B(j)-W(j))_{j=0..n+1}$ are given under both
  configurations and graphically represented.}
\label{fig:number}
\end{figure}

\begin{proof} Let $\Gamma_{n+1}$ be the set of (unconstrained)
  configurations of $n+1$ black and $n+1$ white particles distributed
  between two rows of $n+1$ cells, so that $|\Gamma_{n+1}|={2n+2
    \choose n+1}$. Among these configurations, we restrict our
  attention to the subset $\overline{\Gamma}_{n+1}$ of those ending
  with $|\bw|$- or a $|\bb|$-column.  Exchanging $\bullet$ and $\circ$
  particles is  a bijection between $\overline{\Gamma}_{n+1}$
  and its complement in $\Gamma_{n+1}$, so that
  $|\overline{\Gamma}_{n+1}|=\frac{1}{2}{2n+2\choose n+1}$.
  
  The proof of the lemma consists in a bijection $\phi$ between
  $\Omega_n$ and $\overline{\Gamma}_{n+1}$ (see
  Figure~\ref{fig:number}). Given $\omega \in \Omega_n$, its image
  $\phi(\omega)$ is obtained as follows: First, if {the number of
    $|\xx|$-columns of $\omega$ is even,} add a $|\bw|$-column at the
  end of $\omega$, otherwise add to it an $|\xx|$-column.  Then
  replace the first half of the $|\xx|$-columns by $|\ww|$-columns,
  and the remaining half by $|\bb|$-columns (from left to right).  By
  construction the resulting $\phi(\omega)$ belongs to
  $\overline{\Gamma}_{n+1}$.  Conversely, consider
  $\gamma\in\overline\Gamma_{n+1}$, and let $d=\min(E(j))$ be the
  \emph{depth} of $\gamma$.  Then set $j_i=\min\{j\mid E(j)=-2i\}$,
  and $j'_i=\max\{j\mid E(j-1)=-2i\}$, for $i=1,\ldots,|d/2|$, and
  define the application $\psi$ that first changes columns $j_i$ and
  $j'_i$ into $|\xx|$-columns for all $i=1,\ldots,|d/2|$, and then
  removes the last column.  By construction the blocks between two of
  the modified columns of $\gamma$ satisfy the positivity condition,
  so that $\psi(\gamma)\in\Omega_{n+1}$. Finally the applications
  $\phi$ and $\psi$ are clearly inverse of each other.
\end{proof}

\smallskip
\noindent\textbf{Lemmas~\ref{lemma:ciclico} 
and~\ref{lem:count-all-refined}.} {\em For any 
  $k+\ell+m=n$, the cardinality of the set $\Omega^\ell_{k,m}$ of
  complete configurations with with $\ell$ $|\xx|$-columns, $k$ black
  and $m$ white particles in the top row is
  $\frac{\ell+1}{n+1}{n+1\choose k}{n+1\choose m}$.  }\par\smallskip

\begin{figure}
\begin{center}
\includegraphics[width=120mm]{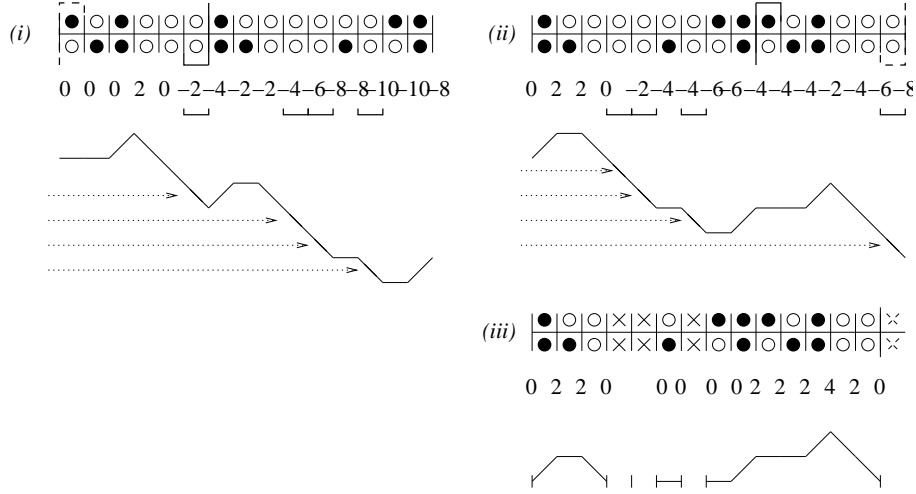}
\end{center}
\caption{(i) An element of $\bar\Delta^{\ell+1}_{k,m}$ (with $\ell=3$
  and column $j=6$ colored), (ii) its conjugate (with column $n+1-j$
  colored), and (iii) the corresponding element of
  $\Omega^\ell_{k,m}$.  The sequence $(B(j)-W(j))_{j=0..n+1}$ is
  given under each configuration and graphically represented.}
\label{fig:ciclico}
\end{figure}

\begin{proof}  The statement is verified using the cycle lemma (see
  \cite[Ch. 11]{lothaire}, or \cite[Ch. 5]{stanley}).  Denote by
  $\Delta^{\ell+1}_{n}$ the set of configurations with $p=n-\ell=k+m$
  black and $p+2\ell+2$ white particles distributed between two rows
  of $n+1$ cells. Then the cardinality of the subset
  $\Delta^{\ell+1}_{k,m}$ of elements of $\Delta^{\ell+1}_{n}$ that
  have $k$ black particles in the top row and the other $m$ in the
  bottom row is $\smash{{n+1 \choose k}{n+1 \choose m}}$. In such a
  configuration the number of white particles exceeds by $2\ell+2$
  that of black particles, so that $E(n+1)=-2\ell-2$.  Given $\omega$
  in $\Delta^{\ell+1}_{k,m}$, let $d=\min(E(j))$ be the depth of
  $\omega$, and set $j_i=\min\{j\mid E(j)=d+2i\}$, for
  $i=0,\ldots,\ell$.  By construction, these $\ell+1$ columns are
  $|\ww|$-columns.  On the one hand, let $\bar{\Delta}^{\ell+1}_{k,m}$
  be the set of pairs $(\omega,j)$ where
  $\omega\in\Delta^{\ell+1}_{k,m}$ and $j\in\{j_0,\ldots,j_\ell\}$, so
  that $|\bar{\Delta}^{\ell+1}_{k,m}|= {n+1 \choose k}{n+1\choose
    m}\cdot (\ell+1)$.  On the other hand, define the set
  $\bar{\Omega}^{\ell+1}_{k,m}$ of pairs $(\omega',i)$ where $\omega'$
  is obtained from an element of $\Omega^{\ell}_{k,m}$ by adding a
  final $|\xx|$-column, and $i\in\{0,\ldots,n\}$. By construction,
  $|\bar{\Omega}^{\ell+1}_{k,m}|= |\Omega^\ell_{k,m}|\cdot (n+1)$.
  
  The proof of the lemma consists in a bijection $\phi$ between
  $\bar{\Delta}^{\ell+1}_{k,m}$ and $\bar{\Omega}^{\ell+1}_{k,m}$ (see
  Figure~\ref{fig:ciclico}). Given $(\omega,j) \in
  \bar{\Delta}^{\ell+1}_{k,m}$, let $\omega_1$ denote the first $j$
  columns of $\omega$, and $\omega_2$ the $n+1-j$ others. Then by
  construction of $j$, the concatenation $\omega_2|\omega_1$ satisfies
  $E(i)>-2\ell-2$ for $i=1,\ldots,n$, and $E(n+1)=-2\ell-2$.  This
  implies that $\omega_2|\omega_1$ decomposes as a sequence
  $\omega'_0,\omega'_1,\ldots,\omega'_{\ell}$ of $\ell+1$ (possibly
  empty) blocks that satisfy the positivity constraint, each followed
  by a $|\ww|$-column.  Let $\omega'$ be obtained by replacing these
  $\ell+1$ $|\ww|$-columns by $|\xx|$-columns.  Then the map
  $(\omega,j)\to(\omega',n+1-j)$ is a bijection of
  $\bar\Delta^{\ell+1}_{k,m}$ onto $\bar{\Omega}^{\ell+1}_{k,m}$: the
  inverse bijection is readily obtained by first replacing the
  $|\xx|$-columns into $|\ww|$-columns, and then recovering the
  factorization $\omega_2|\omega_1$ from the fact that $\omega_2$ has
  $n+1-j$ columns.
\end{proof}

\smallskip
\noindent\textbf{Lemmas~\ref{lemma:olaf} 
and~\ref{lemma:olaf-le-retour}.} {\em 
    The cardinality of the set $\Omega^\ell_{n}$ of complete
    configurations of $\Omega_n$ that have $\ell$ $|\xx|$-columns is
    $\frac{\ell+1}{n+1} { 2n+2 \choose n-\ell}$.
}\par\smallskip
\begin{proof}  The proof uses the same arguments than the proof of
  Lemma \ref{lem:count-all-refined}. The only difference is that, instead of
  counting elements of $\Delta^{\ell+1}_{k,m}$ with $k$ black
  particles in the top row and $m$ in the bottom row, we count
  elements of $\Delta^{\ell+1}_{n}$, the set of configurations of
  $n-\ell$ black particles and $n+2+\ell$ white particles distributed
  in two rows.  Hence the previous factor
  $|\Delta^{\ell+1}_{k,m}|={n+1\choose k}{n+1\choose m}$ is replaced
  by $|\Delta^{\ell+1}_n|={2n+2 \choose n-\ell}$.\end{proof}

\bigskip
\noindent\textbf{Lemma~\ref{lemma:nici}.} {\em The number
  $|\widehat\Omega^\ell_{k,m}|$ of configurations of $|\widehat\Omega_{n}|$
  having $\ell$ $|\xx|$-columns, $k$ black particles at the top, and
  $m$ at the bottom is ${n \choose k}{n \choose m}$.  }\par\smallskip

\begin{proof}  Recall that $\Delta^\ell_{k,m}$ denotes configurations of
  length $n$ with $k$ black and $m+\ell$ white particles in the top
  row, and $m$ black and $k+\ell$ white particles in the bottom row,
  so that $\smash{|\Delta^\ell_{k,m}|={n \choose k}{n \choose m}}$.
  In order to prove the statement of the lemma we show that
  $\Delta^\ell_{k,m}$ and ${\widehat\Omega^\ell_{k,m}}$ are in bijection.
  Let $\smash{\delta \in \Delta^\ell_{k,m}}$, and consider its depth
  $d=\min(E(i))$ and the $\ell$ columns $j_i=\min\{j\mid E(j)=d+2i\}$,
  $i=0,\ldots, \ell-1$, as in the proof of Lemma~\ref{lemma:olaf}.  By
  definition of these columns, the positivity condition is satisfied
  by each block between two of them. Moreover, by definition of $j_0$
  and $j_{\ell-1}$, the positivity condition is also satisfied by the
  concatenation $\omega_\ell|\omega_0$ of the final block
  $\omega_\ell$ and the initial block $\omega_0$.  Hence transforming
  the columns $j_0,\ldots,j_\ell$ into $|\xx|$-columns, and arranging
  the two rows in a circle by fusing walls $0$ and $n$ at the apex
  yields a configuration $\phi(\delta)$ of $\widehat\Omega^\ell_{k,m}$
  (recall that these configurations are not considered up to
  rotation). Conversely, given $\omega$ in $\widehat\Omega^\ell_{k,m}$, a
  unique element $\delta$ of $\Delta^\ell_{k,m}$ such that
  $\phi(\delta)=\omega$ is obtained by opening at the apex and
  transforming $|\xx|$-columns into $|\ww|$-columns.
 \end{proof}

\section{Conclusions and relations to Brownian excursions}
\label{sec:conclusion}
The starting point of this paper was a ``combinatorial Ansatz'': the
stationary distribution of the two particle TASEP with boundaries can
be expressed in terms of Catalan numbers hence should have a nice
combinatorial interpretation.  In our interpretation, configurations
of the TASEP are completed by a (usually hidden) second row in which
particles go back.  In the most interesting case
$\alpha=\beta=\gamma=1$, the resulting chain has a uniform stationary
distribution so that the probability of a given TASEP configuration
just reflects the diversity of possible rows hidden below it.

We do not claim that our combinatorial interpretation is of any
physical relevance. However, apart from explaining the ``magical''
occurrence of Catalan numbers in the problem, it sheds new light on
the recent results of Derrida \emph{et al.}  \cite{Derrida03brownian}
connecting the TASEP with Brownian excursion.
More precisely, using explicit calculations, Derrida \emph{et al.}
show that the density of black particles in configurations of the two
particle TASEP can be expressed in terms of a pair $(e_t,b_t)$ of
independent processes, a Brownian excursion $e_t$ and a Brownian
motion $b_t$. In our interpretation these two quantities appear at the
discrete level, associated to each complete configuration $\omega$ of
$\Omega^0_n$:
\begin{itemize}
\item The role of the Brownian excursion for $\omega$ is played by the
  halved differences $e(i)=\frac12(B(i)-W(i))$ between the number of
  black and white particles sitting on the left of wall $i$, for
  $i=0,\ldots,n$.  By definition of complete configurations,
  $(e(i))_{i=0,\ldots,n}$ is a discrete excursion, that is,
  $e(0)=e(n)=0$, $e(i)\geq0$ and $|e(i)-e(i-1)|\in\{0,1\}$, for
  $i=0,\ldots, n$.
\item The role of the Brownian motion is played for $\omega$ by the
  differences $b(i)=B_{top}(i)-B_{bot}(i)$ between the number of black
  particles sitting in the top and in the bottom row, on the left of
  wall $i$, for $i=0,\ldots, n$. This quantity
  $(b(i))_{i=0,\ldots,n}$ is a discrete walk, with
  $|b(i)-b(i-1)|\in\{0,1\}$ for $i=0,\ldots,n$.
\end{itemize}
Since $e(i)+b(i)=2B_{top}(i)-i$, the functions $e$ and $b$ allow one to
describe the cumulated number of black particles in the top row of a
complete configuration. Accordingly, the density of black particles in
a given segment $(i,j)$ is
$(B_{top}(j)-B_{top}(i))/(j-i)=\frac12+\frac{e(j)-e(i)}{2(j-i)}
+\frac{b(j)-b(i)}{2(j-i)}$. This is a discrete version of the quantity
considered by Derrida \emph{et al.} in \cite{Derrida03brownian}.

Now the two walks $e(i)$ and $b(i)$ are correlated since one is
stationary when the other is not, and vice-versa:
$|e(i)-e(i-1)|+|b(i)-b(i-1)|=1$. Given $\omega$, let
$I_e=\{\alpha_1<\ldots<\alpha_p\}$ be the set of indices of $|\bb|$- and
$|\ww|$-columns, and $I_b=\{\beta_1<\ldots<\beta_q\}$ the set of indices
of $|\bw|$- and $|\wb|$-columns ($p+q=n$). Then the walk
$e'(i)=e(\alpha_i)-e(\alpha_{i-1})$ is the excursion obtained from $e$
by ignoring stationary steps, and the walk
$b'(i)=b(\beta_i)-b(\beta_{i-1})$ is obtained from $b$ in the same
way.  Conversely given a simple excursion $e'$ of length $p$, a simple
walk $b'$ of length $q$ and a subset $I_e$ of $\{1,\ldots,p+q\}$ of
cardinality $p$, two correlated walks $e$ and $b$, and thus a complete
configuration $\omega$ can be uniquely reconstructed.  The consequence
of this discussion is that the uniform distribution on $\Omega^0_n$
corresponds to the uniform distribution of triples $(I_e,e',b')$ where,
given $I_e$, the processes $e'$ and $b'$ are independent.

A direct computation shows that in the large $n$ limit, with
probability exponentially close to~1, a random configuration $\omega$
is described by a pair $(e',b')$ of walks of roughly equal lengths
$n/2+O(n^{1/2+\varepsilon})$.  In particular, up to multiplicative
constants, the normalized pairs
$(\frac{e'(tn/2)}{n^{1/2}},\frac{b'(tn/2)}{n^{1/2}})$ and
$(\frac{e(tn)}{n^{1/2}},\frac{b(tn)}{n^{1/2}})$ both converge 
to the same pair $(e_t,b_t)$ of independent processes, with $e_t$ a standard
Brownian excursion and $b_t$ a standard Brownian walk.  

Another possible outcome of our approach could be an explicit
construction of a continuum TASEP by taking the limit of the Markov
chain $X$, viewed as a Markov chain on pairs of walks.  An appealing
way to give a geometric meaning to the transitions in the continuum
limit could be to use a representation in terms of parallelogram
polyominos \cite{stanley}, using the process $e(t)$ (or $e_t$ in the continuum limit)
to describe the width of the polymonino and the process $b(t)$ (or
$b_t$ in the continuum limit) to describe the vertical displacement of
its spine.

\subsection*{Acknowledgments.} Referees are warmly thanked for their 
great help in improving the paper.

\bibliography{libre}

%%%%%%%%%%%%%%%%%%%%%%%%%%%

\end{document}